\def\BibCommonsStandalone{}

\documentclass[11pt]{article}
\usepackage{adjustbox,aliascnt,amsmath,amssymb,amsthm,calc,colonequals,enumitem,mathtools,mleftright,parskip,xcolor}
\usepackage[a4paper,left=1.3in,right=1.3in,top=1.25in,bottom=1.25in,footskip=0.6in]{geometry}
\usepackage[unicode]{hyperref}
\usepackage[en-GB]{datetime2}

\usepackage{bib-commons}

\newcommand{\newaliastheorem}[2]{%
  \newaliascnt{#1}{theorem}
  \newtheorem{#1}[#1]{#2}
  \aliascntresetthe{#1}
  \expandafter\def\csname #1autorefname\endcsname{#2}
}

\theoremstyle{plain}
\newtheorem{theorem}{Theorem}[section]
\newtheorem*{theorem*}{Theorem}
\newaliastheorem{lemma}{Lemma}
\newtheorem*{lemma*}{Lemma}
\newaliastheorem{proposition}{Proposition}
\newtheorem*{proposition*}{Proposition}
\newaliastheorem{corollary}{Corollary}
\newtheorem*{corollary*}{Corollary}
\newtheorem*{fact*}{Fact}

\theoremstyle{definition}
\newtheorem*{definition}{Definition}

\newcommand{\tuple}[1]{\left\langle #1 \right\rangle}

\newcommand{\Seq}{\mathrm{Seq}}
\newcommand{\Rec}{\mathrm{Rec}}

\DeclareMathOperator{\dom}{dom}
\DeclareMathOperator{\ran}{ran}
\DeclareMathOperator{\field}{field}

\DeclareMathOperator{\otyp}{otyp}

\DeclareMathOperator{\TJ}{TJ}

\DeclareMathOperator{\WO}{WO}
\DeclareMathOperator{\KB}{KB}

\title{An ordinal analysis of $\mathrm{CM}$ and its extensions}
\author{Shuwei Wang}
\date{\DTMdate{2025-01-21}}

\begin{document}

\maketitle

\begin{abstract}
    In \cite{weaver09-cm}, Nik Weaver proposed a novel intuitionistic formal theory of third-order arithmetic as a formalisation of his philosophical position known as \emph{mathematical conceptualism}. In this paper, we will construct a realisability model from the partial combinatory algebra of $\Sigma^1_1$-definable partial functions and use it to provide an ordinal analysis of this formal theory. Additionally, we will examine possible extensions to this system by adding well-ordering axioms, which are briefly mentioned but never thoroughly studied in Weaver's work. We aim to use the realisability arguments to discuss how much such extensions constitute an increase from the original theory's proof-theoretic strength.
\end{abstract}

\tableofcontents

\section{Introduction}

Nik Weaver first introduced his philosophy of \emph{mathematical conceptualism} in \cite{weaver05-conceptualism}, which was essentially rejecting the full Cantorian universe of sets as a foundational ontological commitment and instead advocating a restricted language involving only countably and predicatively constructed mathematical objects --- for example, Weaver has no problem reasoning about individual real numbers, but views the complete totality $\mathcal{P}\mleft(\mathbb{N}\mright)$ as a dangerous object to postulate.

In his early papers \cites{weaver05-conceptualism,weaver05-j2}, this philosophical motivation is only reflected by a mild restriction of set theory to the weak fragment $J_2$, the first two levels in Jensen's hierarchy in the fine structure theory, but Weaver later adopted a more radical foundational viewpoint that quantifiers over infinitary objects (e.g.\ the individual reals) can be manipulated only in an intuitionistic manner. This led to a shift from using a set-theoretic foundation to one in higher order arithmetic, which Weaver formalised in \cite{weaver09-cm}. The system consists of the following axioms over third-order intuitionistic logic:
\begin{enumerate}[label=$\left(\arabic*\right)$]
    \item the number axioms, essentially the axioms of $\mathrm{PA}$ excluding induction;
    \item the unrestricted induction scheme
          \[\varphi\mleft(0\mright) \land \forall n \left(\varphi\mleft(n\mright) \rightarrow \varphi\mleft(n + 1\mright)\right) \rightarrow \forall n \, \varphi\mleft(n\mright);\]
          for any third-order formula $\varphi\mleft(n\mright)$;
    \item the recursion scheme, also known as dependent choice,
          \[\forall n \, \forall X \, \exists Y \, \varphi\mleft(n, X, Y\mright) \rightarrow \forall X \, \exists Z \left(\left(Z\right)_0 = X \land \forall n \, \varphi\mleft(n, \left(Z\right)_n, \left(Z\right)_{n + 1}\mright)\right),\]
          for any third-order formula $\varphi\mleft(n, X, Y\mright)$, where uppercase $X, Y, Z$ denote second-order variables, and $\left(Z\right)_n$ denotes the set $\left\{m : \tuple{n, m} \in Z\right\}$ (which exists by the comprehension axioms below);
    \item the decidable comprehension axioms
          \begin{align*}
              \forall n \left(\varphi\mleft(n\mright) \lor \neg \varphi\mleft(n\mright)\right) & \rightarrow \exists X \, \forall n \left(n \in X \leftrightarrow \varphi\mleft(n\mright)\right),                   \\
              \forall X \left(\varphi\mleft(X\mright) \lor \neg \varphi\mleft(X\mright)\right) & \rightarrow \exists \mathbf{X} \, \forall X \left(X \in \mathbf{X} \leftrightarrow \varphi\mleft(X\mright)\right),
          \end{align*}
          where $X$ denotes a second-order variable and boldface $\mathbf{X}$ denotes a third-order variable;
    \item and lastly, the law of excluded middle $\varphi \lor \neg \varphi$ for all arithmetic formulae $\varphi$, i.e.\ formulae that do not contain second- or third-order quantifiers.
\end{enumerate}

Our paper will be more concerned with Weaver's formal theory itself than with the philosophical justifications behind it. It is already observed by Rathjen in \cite{rathjen23-cst} that his theory is interpretable in the weak set theory $\mathrm{CZF} + \mathrm{LPO} + \mathrm{RDC}$, which is known to have the same proof-theoretic ordinal as $\mathrm{CZF}$ itself, namely the Bachmann--Howard ordinal. However, we will demonstrate a much lower upper bound for its strength in \autoref{sec:cm-ordinal-analysis}, where a realisability argument is constructed to show that

\newcommand{\AC}{\text{$\Sigma^1_1$-$\mathrm{AC}$}}

\begin{theorem}
    The classical fragment $\AC$ of second-order arithmetic (with full induction) interprets Weaver's system $\mathrm{CM}$.
\end{theorem}

This implies that the proof-theoretic ordinal of $\mathrm{CM}$ is precisely $\varphi_{\varepsilon_0}\mleft(0\mright)$.

In \autoref{sec:axiom-gwo}, we will look at the implications of extending Weaver's system by asserting that all second-order objects form a global well-ordering, which Weaver proposed originally in section 2.3 of \cite{weaver09-cm} but did not discuss in detail afterwards. The rest of the paper will then be concerned with determining the effect of possible extensions on the proof-theoretic strength of the theory: in \autoref{sec:cm-gwo-ordinal-analysis}, we will show that the classical theory $\mathrm{BI}$ of bar induction is required to interpret a full axiomatisation of the global well-ordering, rendering its proof-theoretic ordinal as the Bachmann--Howard ordinal $\theta_{\varepsilon_{\Omega + 1}}\mleft(0\mright)$; while in \autoref{sec:countable-initial-segments}, we will additionally show that if the axiom asserting the countability of the initial segments of this well-ordering is omitted, the same strength cannot be achieved.

\section{An ordinal analysis of \texorpdfstring{$\mathrm{CM}$}{CM}}
\label{sec:cm-ordinal-analysis}

\subsection{The axioms of \texorpdfstring{$\mathrm{CM}$}{CM}}
\label{subsec:cm-axioms}

To begin, we will first make a few clarifications/modifications to Weaver's formal theory, for the convenience of our later analysis. Following Weaver, we still treat the theory $\mathrm{CM}$ in a third-order formal language where first-order variables are denoted by lowercase letters $m, n, \ldots$, second-order variables are denoted by italic uppercase letters $X, Y, \ldots$, and third-order variables are denoted by bold uppercase letters $\mathbf{X}, \mathbf{Y}, \ldots$. The membership relation between first-and-second-order and second-and-third-order objects will be distinguished\footnote{This is done so that, when we work with a realisability model and have natural number realisers for all three orders of variables, we can still read off the order of the realiser from the presentation of the formula directly.} and rendered as $n \in_1 X$ and $X \in_2 \mathbf{X}$.

Non-logical symbols of this language shall consist of all primitive symbols from the language of $\mathrm{PA}$ (i.e.\ zero, successors, addition and multiplication) as well as a binary relation symbol $\prec$ between two second-order variables. This relation will serve as the global well-ordering on second-order objects in the extensions, so for convenience we will enforce that it is present in the formal language from the beginning, though the only axiom in $\mathrm{CM}$ concerning it will be to assert that this is a decidable relation:
\begin{enumerate}
    \item[$\left(6\right)$] $\forall X \, \forall Y \left(X \prec Y \lor \neg X \prec Y\right)$.
\end{enumerate}

In the language, we will treat first-order equality as a primitive symbol with the standard equality axioms, while second- and third-order equalities are treated as abbreviations:
\begin{align*}
    X_1 = X_2 \ \                    & \ratio\Leftrightarrow \ \ \forall n \left(n \in_1 X_1 \leftrightarrow n \in_1 X_2\right),                   \\
    \mathbf{X}_1 = \mathbf{X}_2 \ \  & \ratio\Leftrightarrow \ \ \forall X \left(X \in_2 \mathbf{X}_1 \leftrightarrow X \in_2 \mathbf{X}_2\right).
\end{align*}

Although not mentioned in Weaver's formalisation, we will include the following axioms asserting the extensionality of relations and predicates for second-order variables:
\begin{enumerate}
    \item[$\left(7\right)$] $X = Y \rightarrow \left(X \in_2 \mathbf{X} \leftrightarrow X \in_2 \mathbf{Y}\right)$,
    \item[$\left(8\right)$] $X = Y \land X' = Y' \rightarrow \left(X \prec Y \leftrightarrow X' \prec Y'\right)$.
\end{enumerate}
These are needed in order to show that

\begin{proposition}[]
    Given axioms $\left(7\right)$ and $\left(8\right)$, the Leibniz law holds, i.e.
    \[t = s \rightarrow \left(\varphi\mleft(t\mright) \leftrightarrow \varphi\mleft(s\mright)\right),\]
    for any formula $\varphi$ and any terms $t, s$ of the same, correct order.
\end{proposition}

\begin{proof}
    This follows from simple induction on the complexity of $\varphi$. Just notice that all atomic cases hold trivially by definition of the equality symbol (which is primitive for first-order terms, but non-primitive for higher-order terms), except when $\varphi\mleft(t\mright)$ is $t \in_2 \mathbf{X}$, $t \prec X$ or $X \prec t$, where the extensionality assertions above is used.
\end{proof}

The Leibniz law is necessary for expressions with set-comprehension terms, such as $\left\{n : \varphi\mleft(n\mright)\right\} \in_2 \mathbf{X}$, to act as non-ambiguous abbreviations in the language.

Additionally, the arithmetic decidability scheme $\left(5\right)$ can be complicated to work with at times, so we will break it down into two parts:
\begin{enumerate}
    \item[$\left(5'\right)$] atomic decidability axioms
          \[n = m \lor \neg n = m, \qquad n \in_1 X \lor \neg n \in_1 X, \qquad X \in_2 \mathbf{X} \lor \neg X \in_2 \mathbf{X};\]
    \item[$\left(5''\right)$] the numerical omniscience scheme, also known as the limited principle of omniscience,
          \[\forall n \left(\varphi\mleft(n\mright) \lor \psi\mleft(n\mright)\right) \rightarrow \left(\forall n \, \varphi\mleft(n\mright) \lor \exists n \, \psi\mleft(n\mright)\right)\]
          for any third-order formulae $\varphi\mleft(n\mright)$ and $\psi\mleft(n\mright)$.
\end{enumerate}

We know that

\begin{proposition}
    Over the rest of the axioms of $\mathrm{CM}$, the axiom scheme $\left(5\right)$ is equivalent to $\left(5'\right) + \left(5''\right)$.
\end{proposition}

\begin{proof}
    The forward direction is proved by Weaver as theorem 2.2 in \cite{weaver09-cm}; the backward direction is trivial through an induction on the complexity of an arithmetic formula.
\end{proof}

Therefore, in this paper, we will take $\mathrm{CM}$ to formally denote the collection of axioms $\left(1\right) + \left(2\right) + \left(3\right) + \left(4\right) + \left(5'\right) + \left(5''\right) + \left(6\right) + \left(7\right) + \left(8\right)$, each precisely as defined above.

Lastly, the binary pairing function $\tuple{{-}, {-}} : \mathbb{N} \times \mathbb{N} \rightarrow \mathbb{N}$ will always denote a fixed bijection in this paper, e.g.\ given by
\[\tuple{m, n} = \left(m + n\right)\left(m + n + 1\right) / 2 + m.\]
with inverses $\left(-\right)_0, \left(-\right)_1 : \mathbb{N} \rightarrow \mathbb{N}$. The fact that every natural number is a code for some pair will be convenient for some constructions below. In \cite{weaver09-cm}, Weaver notably used a diffferent, non-surjective map $\tuple{m, n} = 2^m \cdot 3^n$. This difference does affect the precise meaning of certain formulae, for example the axiom scheme $\left(3\right)$ of dependent choice. However, since $\mathrm{CM}$ contains classical first-order $\mathrm{PA}$ as a subtheory, it can be shown that different interpretations of the pairing function are ultimately equivalent.

\subsection{\texorpdfstring{$\Sigma^1_1$}{Sigma(1, 1)} function codes}

We will set up the realisability structure for $\mathrm{CM}$ in a classical context of second-order arithmetic, using $\mathrm{ACA}_0$ as the base theory. Following definition VIII.3.18 in Simpson \cite{simpson09-soa}, we let $\pi\mleft(e, m, a, b, X\mright)$ denote a universal $\Pi^0_1$-formula with all free variables indicated, and correspondingly a universal $\Sigma^1_1$-formula is given by $\exists X \, \pi\mleft(e, m, a, b, X\mright)$, i.e.\ for any $\Sigma^1_1$-formula $\varphi\mleft(m, a, b\mright)$ with all free variables indicated, $\mathrm{ACA}_0$ proves that there exists $e \in \mathbb{N}$ such that
\[\forall m \, \forall a \, \forall b \left(\varphi\mleft(m, a, b\mright) \leftrightarrow \exists X \, \pi\mleft(e, m, a, b, X\mright)\right).\]

We shall use the notation
\[\mleft\{e\mright\}\mleft(a\mright) = b \ \ \ratio\Leftrightarrow \ \ \exists X \, \pi\mleft(e_0, e_1, a, b, X\mright).\]

We say that $\mleft\{e\mright\}\mleft(a\mright)$ \emph{converges}, denoted $\mleft\{e\mright\}\mleft(a\mright) {\downarrow}$, when $\exists!b \, \mleft\{e\mright\}\mleft(a\mright) = b$ holds. We write
\begin{align*}
    \dom\mleft(e\mright)    & \coloneqq \left\{a : \mleft\{e\mright\}\mleft(a\mright) {\downarrow}\right\},                           \\
    \ran\mleft(e, X\mright) & \coloneqq \left\{b : \exists a \in X \, \mleft\{e\mright\}\mleft(a\mright) = b\right\},                 \\
    \left[X, Y\right]       & \coloneqq \left\{e : X \subseteq \dom\mleft(e\mright) \land \ran\mleft(e, X\mright) \subseteq Y\right\}
\end{align*}
for any definable classes $X, Y \subseteq \mathbb{N}$. Note that $\dom\mleft(e\mright)$ is $\Sigma^1_2$ and hence so is $\left[X, Y\right]$, even when $X$ and $Y$ are both arithmetic subsets of $\mathbb{N}$.

It follows immediately that
\[\varphi_0\mleft(m, a, b\mright) \ \ \ratio\Leftrightarrow \ \ \mleft\{a_0\mright\}\mleft(a_1\mright) = b\]
is a (parameterless) $\Sigma^1_1$-formula and correspondingly, we can find $e \in \mathbb{N}$ so that for any $f \in \mathbb{N}$ and $a \in \dom\mleft(f\mright)$, we have $\tuple{f, a} \in \dom\mleft(\tuple{e, 0}\mright)$ and $\mleft\{\tuple{e, 0}\mright\}\mleft(\tuple{f, a}\mright) = \mleft\{f\mright\}\mleft(a\mright)$.

In fact, we can prove the following parametrised generalisation:

\begin{proposition}
    \label{prop:sigma11-smn}
    Given any definable class $A$ of natural numbers and a $\Sigma^1_1$-formula $\varphi\mleft(m, a, b\mright)$ such that, for any $a \in A$, there exists a unique $b \in \mathbb{N}$ satisfying $\varphi\mleft(a_0, a_1, b\mright)$, then there exists $e$ such that for any $a \in A$, we have $a_0 \in \dom\mleft(e\mright)$, $a_1 \in \dom\mleft(\mleft\{e\mright\}\mleft(a_0\mright)\mright)$ and
    \[\varphi\mleft(a_0, a_1, \mleft\{\mleft\{e\mright\}\mleft(a_0\mright)\mright\}\mleft(a_1\mright)\mright).\]
\end{proposition}

\begin{proof}
    Simply take $f, n \in \mathbb{N}$ such that
    \begin{align*}
        \forall m \, \forall a \, \forall b \left(b = \tuple{m, a} \leftrightarrow \exists X \, \pi\mleft(f, m, a, b, X\mright)\right), \\
        \forall m \, \forall a \, \forall b \left(\varphi\mleft(m, a, b\mright) \leftrightarrow \exists X \, \pi\mleft(n, m, a, b, X\mright)\right).
    \end{align*}
    If follows that $\mleft\{\mleft\{\tuple{f, n}\mright\}\mleft(m\mright)\mright\}\mleft(a\mright) = b$ iff $\mleft\{\tuple{n, m}\mright\}\mleft(a\mright) = b$ iff $\varphi\mleft(m, a, b\mright)$.
\end{proof}

\begin{corollary}
    \label{cor:sigma11-pca}
    The set of natural numbers $\mathbb{N}$ form a partial combinatory algebra (with arithmetic)\footnote{By a \emph{partial combinatory algebra with arithmetic}, we require that the distinguished constants $k, s, 0, \mathrm{succ}, \mathrm{pred}, d, p, p_0, p_1$ are defined in the structure for, respectively, the two basic combinators, zero, successor and predecessor on $\mathbb{N}$, definition by cases and the pairing and unpairing functions. Such PCAs are first considered by Feferman in \cite{feferman75-em} as part of a stronger system for explicit mathematics.} under the application operation $\mleft\{-\mright\}\mleft(-\mright)$.
\end{corollary}

\begin{proof}
    The only non-trivial tasks are to define the basic combinators $k$ and $s$. For $k$, we directly invoke \autoref{prop:sigma11-smn} to compute some $k \in \mathbb{N}$ such that for every $a \in \mathbb{N}$ and $b \in \dom\mleft(a\mright)$, we have $a \in \dom\mleft(k\mright)$, $b \in \dom\mleft(\mleft\{k\mright\}\mleft(a\mright)\mright)$ and
    \[\mleft\{\mleft\{k\mright\}\mleft(a\mright)\mright\}\mleft(a\mright) = \mleft\{a\mright\}\mleft(b\mright)\]
    holds. For $s$, we first use \autoref{prop:sigma11-smn} to compute $s' \in \mathbb{N}$ such that for every $a, b \in \mathbb{N}$ and $c \in \dom\mleft(a\mright) \cap \dom\mleft(b\mright)$ such that $\mleft\{b\mright\}\mleft(c\mright) \in \dom\mleft(\mleft\{a\mright\}\mleft(c\mright)\mright)$, we have $\tuple{a, b} \in \dom\mleft(s'\mright)$, $c \in \dom\mleft(\mleft\{s'\mright\}\mleft(\tuple{a, b}\mright)\mright)$ and
    \[\mleft\{\mleft\{s'\mright\}\mleft(\tuple{a, b}\mright)\mright\}\mleft(c\mright) = \mleft\{\mleft\{a\mright\}\mleft(c\mright)\mright\}\mleft(\mleft\{b\mright\}\mleft(c\mright)\mright)\]
    holds. Lastly, we use \autoref{prop:sigma11-smn} again to compute $s \in \mathbb{N}$ such that for any $a, b \in \mathbb{N}$ (where at least one aforementioned $c$ exists), we have $a \in \dom\mleft(s\mright)$, $b \in \dom\mleft(\mleft\{s\mright\}\mleft(a\mright)\mright)$ and
    \[\mleft\{\mleft\{s\mright\}\mleft(a\mright)\mright\}\mleft(b\mright) = \mleft\{s'\mright\}\mleft(\tuple{a, b}\mright). \qedhere\]
\end{proof}

\subsection{Realisability of \texorpdfstring{$\mathrm{CM}$}{CM}}
\label{subsec:realisability-cm}

We shall show that the axiom system $\mathrm{CM}$, as defined in \autoref{subsec:cm-axioms}, can be interpreted in $\AC$ (with unrestricted induction). This interpretation will be set-up as a realisability model where assignments for first-order variables range over $\mathbb{N}$, assignments for second-order variables range over $S_2 = \left[\mathbb{N}, \left\{0, 1\right\}\right]$, and assignments for third-order variables range over
\[S_3 = \left\{d \in \left[S_2, \left\{0, 1\right\}\right] : \forall e, f \in S_2 \left(e \equiv f \rightarrow \mleft\{d\mright\}\mleft(e\mright) = \mleft\{d\mright\}\mleft(f\mright)\right)\right\},\]
where the binary relation $\equiv$ on $S_2$ is defined as
\[e \equiv f \ \ \ratio\Leftrightarrow \ \ \forall n \, \mleft\{e\mright\}\mleft(n\mright) = \mleft\{f\mright\}\mleft(n\mright).\]
Fix an arbitrary interpretation $r \in \left[S_2 \times S_2, \left\{0, 1\right\}\right]$ for the relation $\prec$, also satisfying
\[\forall e, e', f, f' \in S_2 \left(e \equiv f \land e' \equiv f' \rightarrow \mleft\{r\mright\}\mleft(\tuple{e, e'}\mright) = \mleft\{r\mright\}\mleft(\tuple{f, f'}\mright)\right),\]
for example just $r = \lambda n. 0$, we shall recursively define the following realisability statements:
\begin{align*}
    d \Vdash t = s \ \                                                   & \ratio\Leftrightarrow \ \ t = s \text{ for any arithmetic terms $t$ and $s$},                                                                                  \\
    d \Vdash t \in_1 e \ \                                               & \ratio\Leftrightarrow \ \ \mleft\{e\mright\}\mleft(t\mright) = 1 \text{ for any arithmetic terms $t$},                                                         \\
    d \Vdash e \in_2 f \ \                                               & \ratio\Leftrightarrow \ \ \mleft\{f\mright\}\mleft(e\mright) = 1,                                                                                              \\
    d \Vdash e \prec f \ \                                               & \ratio\Leftrightarrow \ \ \mleft\{r\mright\}\mleft(\tuple{e, f}\mright) = 1,                                                                                   \\
    d \Vdash \varphi \land \psi \ \                                      & \ratio\Leftrightarrow \ \ d_0 \Vdash \varphi \land d_1 \Vdash \psi,                                                                                            \\
    d \Vdash \varphi \lor \psi \ \                                       & \ratio\Leftrightarrow \ \ \left(d_0 = 0 \land d_1 \Vdash \varphi\right) \lor \left(d_0 = 1 \land d_1 \Vdash \psi\right),                                       \\
    d \Vdash \neg \varphi \ \                                            & \ratio\Leftrightarrow \ \ \forall e \, \neg e \Vdash \varphi,                                                                                                  \\
    d \Vdash \varphi \rightarrow \psi \ \                                & \ratio\Leftrightarrow \ \ \forall e \left(e \Vdash \varphi \rightarrow e \in \dom\mleft(d\mright) \land \mleft\{d\mright\}\mleft(e\mright) \Vdash \psi\right), \\
    d \Vdash \forall x \, \varphi\mleft(x\mright) \ \                    & \ratio\Leftrightarrow \ \ d \in \left[\mathbb{N}, \mathbb{N}\right] \land \forall n \, \mleft\{d\mright\}\mleft(n\mright) \Vdash \varphi\mleft(n\mright),      \\
    d \Vdash \exists x \, \varphi\mleft(x\mright) \ \                    & \ratio\Leftrightarrow \ \ d_1 \Vdash \varphi\mleft(d_0\mright),                                                                                                \\
    d \Vdash \forall X \, \varphi\mleft(X\mright) \ \                    & \ratio\Leftrightarrow \ \ d \in \left[S_2, \mathbb{N}\right] \land \forall n \in S_2 \, \mleft\{d\mright\}\mleft(n\mright) \Vdash \varphi\mleft(n\mright),     \\
    d \Vdash \exists X \, \varphi\mleft(X\mright) \ \                    & \ratio\Leftrightarrow \ \ d_0 \in S_2 \land d_1 \Vdash \varphi\mleft(d_0\mright),                                                                              \\
    d \Vdash \forall \mathbf{X} \, \varphi\mleft(\mathbf{X}\mright) \ \  & \ratio\Leftrightarrow \ \ d \in \left[S_3, \mathbb{N}\right] \land \forall n \in S_3 \, \mleft\{d\mright\}\mleft(n\mright) \Vdash \varphi\mleft(n\mright),     \\
    d \Vdash \exists \mathbf{X} \, \varphi\mleft(\mathbf{X}\mright) \ \  & \ratio\Leftrightarrow \ \ d_0 \in S_3 \land d_1 \Vdash \varphi\mleft(d_0\mright).
\end{align*}

For any third-order formula $\varphi\mleft(x_1, \ldots, x_n, X_1, \ldots, X_m, \mathbf{X}_1, \ldots, \mathbf{X}_r\mright)$ with all free variables indicated, we write $d \Vdash \varphi$ to denote
\begin{align*}
     & \forall a_1, \ldots, a_n \in \mathbb{N} \, \forall b_1, \ldots, b_m \in S_2 \, \forall c_1, \ldots, c_r \in S_3                                                                         \\
     & \qquad \mleft\{d\mright\}\mleft(\tuple{a_1, \ldots, a_n, b_1, \ldots, b_m, c_1, \ldots, c_r}\mright) \Vdash \varphi\mleft(a_1, \ldots, a_n, b_1, \ldots, b_m, c_1, \ldots, c_r\mright).
\end{align*}

We shall now show in $\AC$ that whenever
\[\mathrm{CM} \vdash \varphi\mleft(x_1, \ldots, x_n, X_1, \ldots, X_m, \mathbf{X}_1, \ldots, \mathbf{X}_r\mright),\]
$\varphi$ is realised in the sense above.

Firstly, it is easy to verify that our realisability conditions interpret intuitionistic logic. Take note of the following before we begin:
\begin{itemize}
    \item Since we work in $\AC$, the pointclass $\Sigma^1_1$ is closed under adding any first-order quantifications.
    \item Here we treat third-order logic simply as a $3$-sorted language in first-order logic, therefore all logical axioms and rules of inferences have to be verified for all three sorts. Luckily, due to the way we defined the realisability conditions above, verification for all the sorts are entirely analogous.
\end{itemize}

Now, we use the same Hilbert-style formal system as Weaver mentioned in \cite{weaver09-cm}. By \autoref{cor:sigma11-pca}, the realisers can simply be given as $\lambda$-constructs:
\begin{itemize}
    \item The axiom $\varphi \rightarrow \left(\psi \rightarrow \varphi\right)$ is realised\footnote{Here and in all subsequent realisability claims, the first parameter $p$ accommodates for the tuple of assignments $\tuple{a_1, \ldots, a_n, b_1, \ldots, b_m, c_1, \ldots, c_r}$ to free variables in the axiom.} by $\lambda p. \lambda e. \lambda f. e$;
    \item $\left(\varphi \rightarrow \psi\right) \rightarrow \left(\left(\varphi \rightarrow \left(\psi \rightarrow \sigma\right)\right) \rightarrow \left(\varphi \rightarrow \sigma\right)\right)$ is realised by
          \[\lambda p. \lambda f. \lambda g. \lambda e. \mleft\{\mleft\{g\mright\}\mleft(e\mright)\mright\}\mleft(\mleft\{f\mright\}\mleft(e\mright)\mright);\]
    \item $\varphi \rightarrow \left(\psi \rightarrow \varphi \land \psi\right)$ is realised by $\lambda p. \lambda e. \lambda f. \tuple{e, f}$;
    \item $\varphi \land \psi \rightarrow \varphi$ is realised by $\lambda p. \lambda e. e_0$;
    \item $\varphi \land \psi \rightarrow \psi$ is realised by $\lambda p. \lambda e. e_1$.
    \item $\varphi \rightarrow \varphi \lor \psi$ is realised by $\lambda p. \lambda e. \tuple{0, e}$;
    \item $\psi \rightarrow \varphi \lor \psi$ is realised by $\lambda p. \lambda e. \tuple{1, e}$;
    \item $\left(\varphi \rightarrow \sigma\right) \rightarrow \left(\left(\psi \rightarrow \sigma\right) \rightarrow \left(\varphi \lor \psi \rightarrow \sigma\right)\right)$ is realised by
          \[\lambda p. \lambda f. \lambda g. \lambda n. \left\{\begin{aligned}
                   & \mleft\{f\mright\}\mleft(n_1\mright) &  & \text{if} \ n_0 = 0, \\
                   & \mleft\{g\mright\}\mleft(n_1\mright) &  & \text{otherwise};    \\
              \end{aligned}\right.\]
    \item $\left(\varphi \rightarrow \psi\right) \rightarrow \left(\left(\varphi \rightarrow \neg \psi\right) \rightarrow \neg \varphi\right)$ is realised simply by $\lambda p. \lambda f. \lambda g. 0$. This is a correct realiser because, suppose that $f \Vdash \varphi \rightarrow \psi$, $g \Vdash \varphi \rightarrow \neg \psi$, $d \Vdash \varphi$, then $\mleft\{g\mright\}\mleft(d\mright) \Vdash \neg \psi$ while $\mleft\{f\mright\}\mleft(d\mright) \Vdash \psi$, a contradiction. Therefore, no such realiser $d$ exists and $0$ --- in fact, anything --- realises $\neg \varphi$.
    \item $\varphi\mleft(t\mright) \rightarrow \exists x \, \varphi\mleft(x\mright)$ is realised\footnote{Here $t$ can be a term in any of the three sorts, and (if $t$ is a first-order term) it may contain free parameters. The same applies to cases below.} by $\lambda p. \lambda e. \tuple{t, e}$, where parameters in $p$ may be used when $t$ is a complex first-order term;
    \item $\forall x \, \varphi\mleft(x\mright) \rightarrow \varphi\mleft(t\mright)$ is realised by $\lambda p. \lambda f. \mleft\{f\mright\}\mleft(t\mright)$, where again parameters in $p$ may be used when computing $t$.
    \item Regarding the equality axioms for the primitive first-order equality symbol, it suffices to realise $t = t$ by $\lambda p. 0$ and $t = s \land \varphi\mleft(t\mright) \rightarrow \varphi\mleft(s\mright)$ by $\lambda p. \lambda e. e_1$. The latter works since it is easy to show by induction that, when $t$ and $s$ are actually equal terms in $\mathrm{PA}$, $d \Vdash \varphi\mleft(t\mright)$ implies $d \Vdash \varphi\mleft(s\mright)$ as well.
\end{itemize}

For the rules of inferences:
\begin{itemize}
    \item[$\mathbf{MP}$] Given $e \Vdash \varphi$ and $f \Vdash \varphi \rightarrow \psi$, the conclusion $\psi$ is realised by $\lambda p. \mleft\{\mleft\{f\mright\}\mleft(p\mright)\mright\}\mleft(\mleft\{e\mright\}\mleft(p\mright)\mright)$.
    \item[$\forall$-$\mathbf{Gen}$] Given $f \Vdash \varphi \rightarrow \psi$, where $x$ does not occur in $\varphi$, the conclusion $\varphi \rightarrow \forall x \, \psi$ is realised by
          \[\lambda p. \lambda d. \lambda n. \mleft\{\mleft\{f\mright\}\mleft(p, n\mright)\mright\}\mleft(d\mright),\]
          where by $\mleft\{f\mright\}\mleft(p, n\mright)$ we mean that $n$ is passed as an assignment to the free variable $x$ in $\varphi \rightarrow \psi$.
    \item[$\exists$-$\mathbf{Gen}$] Given $f \Vdash \varphi \rightarrow \psi$, where $x$ does not occur in $\psi$, the conclusion $\exists x \, \varphi \rightarrow \psi$ is similarly realised by
          \[\lambda p. \lambda d. \mleft\{\mleft\{f\mright\}\mleft(p, d_0\mright)\mright\}\mleft(d_1\mright).\]
\end{itemize}

It suffices now to further show that the non-logical axioms of $\mathrm{CM}$ are realised. To begin, the number axioms $\left(1\right)$ in $\mathrm{CM}$ are just true, quantifier-free formulae in $\mathrm{PA}$ (with free variables), which can be easily realised by any $\Sigma^1_1$ function with appropriate domain and range.

For the schemes $\left(2\right)$ and $\left(3\right)$, we will need the recursion lemma for $\Sigma^1_1$ functions:

\begin{lemma}[$\AC$]
    There exists $f \in \mathbb{N}$ such that, consider any definable class $A$ of natural numbers and $e \in \left[A, \mathbb{N}\right]$ such that for any $a \in A$, $\mleft\{e\mright\}\mleft(a\mright) \in \left(A\right)_{a_0 + 1}$, then for any $z \in \left(A\right)_0$, we have $\tuple{e, z} \in \dom\mleft(f\mright)$, $d = \mleft\{f\mright\}\mleft(\tuple{e, z}\mright) \in \left[\mathbb{N}, \mathbb{N}\right]$ and
    \[\mleft\{d\mright\}\mleft(0\mright) = z \land \forall n \left(\mleft\{d\mright\}\mleft(n\mright) \in \left(A\right)_n \land \mleft\{d\mright\}\mleft(n + 1\mright) = \mleft\{e\mright\}\mleft(\tuple{n, \mleft\{d\mright\}\mleft(n\mright)}\mright)\right).\]
\end{lemma}

\begin{proof}
    The desirable function computed by $d$ is given by
    \[\eta\mleft(n; e, z\mright) = y \ \ \ratio\Leftrightarrow \ \ \exists s \in \Seq \left(s_0 = z \land s_n = y \land \forall i < n \, \mleft\{e\mright\}\mleft(\tuple{n, s_i}\mright) = s_{i + 1}\right).\]
    where $\Seq$ is the class of all finite sequences coded in $\mathbb{N}$. Therefore, we can find some $e_0 \in \mathbb{N}$ so that $\eta\mleft({-}; e, z\mright)$ is computed by $\Sigma^1_1$ function code $\tuple{e_0, \tuple{e, z}}$, and by \autoref{prop:sigma11-smn}, we can find $f$ to compute $\lambda a. \tuple{e_0, a}$.

    Now, it follows from induction that $f$ does have the desired property. Namely, if $A$ is $\Pi^1_{k + 1}$, then $\Sigma^1_{k + 2}$-induction shows that $\forall n \, \mleft\{\mleft\{f\mright\}\mleft(\tuple{e, z}\mright)\mright\}\mleft(n\mright) \in \left(A\right)_n$.
\end{proof}

\begin{proposition}
    The induction scheme $\left(2\right)$ is realised:
    \[\varphi\mleft(0\mright) \land \forall n \left(\varphi\mleft(n\mright) \rightarrow \varphi\mleft(n + 1\mright)\right) \rightarrow \forall n \, \varphi\mleft(n\mright).\]
\end{proposition}

\begin{proof}
    For any $d \Vdash \varphi\mleft(0\mright)$ and $f \Vdash \forall n \left(\varphi\mleft(n\mright) \rightarrow \varphi\mleft(n + 1\mright)\right)$, a function $g \Vdash \forall n \, \varphi\mleft(n\mright)$ can be given recursively:
    \begin{align*}
        \mleft\{g\mright\}\mleft(0\mright)     & = d,                                                                                                    \\
        \mleft\{g\mright\}\mleft(n + 1\mright) & = \mleft\{\mleft\{f\mright\}\mleft(n\mright)\mright\}\mleft(\mleft\{g\mright\}\mleft(n\mright)\mright).
    \end{align*}
    In this recursion, we will need to check $\mleft\{g\mright\}\mleft(n\mright) \Vdash \varphi\mleft(n\mright)$ for any $n$ and this follows from unrestricted\footnote{Even though we never use an arbitrary amount of second-order quantifiers, the use of unrestrict induction is necessary here: suppose that $\varphi$ contains an alternating sequence of quantifiers, then the statement $d \Vdash \varphi\mleft(p\mright)$ can expands to a form like
        \[d \in \left[\mathbb{N}, \mathbb{N}\right] \land \forall n_0 \, \exists d_1 \in \left[\mathbb{N}, \mathbb{N}\right] \left(\mleft\{d\mright\}\mleft(n_0\mright)_1 = d_1 \land \forall n_1 \, \exists d_2 \in \left[\mathbb{N}, \mathbb{N}\right] \ldots \exists d_i \in \left[\mathbb{N}, \mathbb{N}\right] \ldots \right.\]
        Here the class $\left[\mathbb{N}, \mathbb{N}\right] = \left\{e : \forall a \, \exists! b \, \mleft\{e\mright\}\mleft(a\mright) = b\right\}$ is $\Sigma^1_2$ but lies within a sequence of first-order quantifiers. Since $\AC$ does not show that the pointclass $\Sigma^1_2$ is closed under adding first-order quantifiers, it cannot prove that this statement has bounded complexity.} induction.
\end{proof}

\begin{proposition}
    The recursion scheme $\left(3\right)$ is realised:
    \[\forall n \, \forall X \, \exists Y \, \varphi\mleft(n, X, Y\mright) \rightarrow \forall X \, \exists Z \left(\left(Z\right)_0 = X \land \forall n \, \varphi\mleft(n, \left(Z\right)_n, \left(Z\right)_{n + 1}\mright)\right).\]
\end{proposition}

\begin{proof}
    For any $f \Vdash \forall n \, \forall X \, \exists Y \, \varphi\mleft(n, X, Y\mright)$ and any $x \in S_2$, we first compute some function $g$ recursively:
    \begin{align*}
        \mleft\{g\mright\}\mleft(0\mright)     & = x,                                                                                                      \\
        \mleft\{g\mright\}\mleft(n + 1\mright) & = \mleft\{\mleft\{f\mright\}\mleft(n\mright)\mright\}\mleft(\mleft\{g\mright\}\mleft(n\mright)\mright)_0.
    \end{align*}
    The desired $Z$ is then constructed as some $z \in S_2$ given by $\mleft\{z\mright\}\mleft(\tuple{m, n}\mright) = \mleft\{\mleft\{g\mright\}\mleft(m\mright)\mright\}\mleft(n\mright)$. Notice that $\left(Z\right)_0 = X$ is realised by any function with appropriate domain and range because we are working with implications like $\tuple{0, x} \in_1 Z \rightarrow x \in_1 X$, where the atoms are realised by anything as long as they are true. The conditions $\varphi\mleft(n, \left(Z\right)_n, \left(Z\right)_{n + 1}\mright)$ are realised\footnote{To be completely formal, we may understand $\varphi\mleft(n, \left(Z\right)_n, \left(Z\right)_{n + 1}\mright)$ as
        \[\exists X = \left(Z\right)_n \, \exists Y = \left(Z\right)_{n + 1} \, \varphi\mleft(n, X, Y\mright),\]
        so the formula is realised by $\tuple{\mleft\{g\mright\}\mleft(n\mright), \tuple{*, \tuple{\mleft\{g\mright\}\mleft(n + 1\mright), \tuple{*, \mleft\{\mleft\{f\mright\}\mleft(n\mright)\mright\}\mleft(\mleft\{g\mright\}\mleft(n\mright)\mright)_1}}}}$, where $*$ can be any function code with appropriate domain and range.} by $\mleft\{\mleft\{f\mright\}\mleft(n\mright)\mright\}\mleft(\mleft\{g\mright\}\mleft(n\mright)\mright)_1$.
\end{proof}

The decidable comprehension axioms $\left(4\right)$ are easy cases here:

\begin{proposition}
    The following comprehension axioms are realised:
    \begin{align*}
        \forall n \left(\varphi\mleft(n\mright) \lor \neg \varphi\mleft(n\mright)\right) & \rightarrow \exists X \, \forall n \left(n \in_1 X \leftrightarrow \varphi\mleft(n\mright)\right),                   \\
        \forall X \left(\varphi\mleft(X\mright) \lor \neg \varphi\mleft(X\mright)\right) & \rightarrow \exists \mathbf{X} \, \forall X \left(X \in_2 \mathbf{X} \leftrightarrow \varphi\mleft(X\mright)\right).
    \end{align*}
\end{proposition}

\begin{proof}
    For any $f \Vdash \forall n \left(\varphi\mleft(n\mright) \lor \neg \varphi\mleft(n\mright)\right)$, $X$ can simply be constructed as some $x \in S_2$ such that $\mleft\{x\mright\}\mleft(n\mright) = 1 - \mleft\{f\mright\}\mleft(n\mright)_0$. Observe that $n \in_1 X \rightarrow \varphi\mleft(n\mright)$ is realised by the function $\lambda a. \mleft\{f\mright\}\mleft(n\mright)_1$, where the opposite $\varphi\mleft(n\mright) \rightarrow n \in_1 X$ is realised by an arbitrary function.

    The third-order case is entirely analogous.
\end{proof}

For the atomic decidability axioms $\left(5'\right)$:
\begin{itemize}
    \item $t = s \lor \neg t = s$ is realised by
          \[\lambda p. \left\{\begin{aligned}
                   & \tuple{0, 0} &  & \text{if} \ t = s, \\
                   & \tuple{1, 0} &  & \text{otherwise},
              \end{aligned}\right.\]
          where $p$ represents free parameters used to compute $t$ and $s$;
    \item $t \in_1 e \lor \neg t \in_1 e$ is realised by $\lambda {\tuple{p, e}}. \tuple{1 - \mleft\{e\mright\}\mleft(t\mright), 0}$ where $p$ represents free parameters in $t$;
    \item $e \in_2 f \lor \neg e \in_2 f$ is realised by $\lambda {\tuple{e, f}}. \tuple{1 - \mleft\{f\mright\}\mleft(e\mright), 0}$.
\end{itemize}
And likewise for the decidability of the order relation $\left(6\right)$: $e \prec f \lor \neg e \prec f$ is realised by $\lambda {\tuple{e, f}}. \tuple{1 - \mleft\{r\mright\}\mleft(\tuple{e, f}\mright), 0}$.

\begin{proposition}
    The numerical omniscience scheme $\left(5''\right)$ is realised:
    \[\forall n \left(\varphi\mleft(n\mright) \lor \psi\mleft(n\mright)\right) \rightarrow \left(\forall n \, \varphi\mleft(n\mright) \lor \exists n \, \psi\mleft(n\mright)\right).\]
\end{proposition}

\begin{proof}
    We first fix some $g \in \left[\mathbb{N}, \mathbb{N}\right]$ that computes $\lambda f. \lambda n. \mleft\{f\mright\}\mleft(n\mright)_1$. Now construct the function
    \begin{align*}
        \mleft\{e\mright\}\mleft(f\mright) = d \ \ratio\Leftrightarrow \  & \left(d_0 = 0 \land \forall n \, \mleft\{f\mright\}\mleft(n\mright)_0 = 0 \land d_1 = \mleft\{g\mright\}\mleft(f\mright)\right) \lor                                                                    \\
                                                                          & \left(d_0 = 1 \land \mleft\{f\mright\}\mleft(d_{10}\mright)_0 = 1 \land \forall n < d_{10} \, \mleft\{f\mright\}\mleft(n\mright)_0 = 0 \land d_{11} = \mleft\{f\mright\}\mleft(d_{10}\mright)_1\right).
    \end{align*}
    Observe that since the base theory satisfies classical logic, a unique result $d$ is always well-defined for any $f \Vdash \forall n \left(\varphi\mleft(n\mright) \lor \psi\mleft(n\mright)\right)$. Then $\lambda p. e$ is a realiser for the desired axiom.
\end{proof}

Finally, observe that in our set-up, we have ensured for any $d \in S_3$,
\[\forall e, f \in S_2 \left(e \equiv f \rightarrow \mleft\{d\mright\}\mleft(e\mright) = \mleft\{d\mright\}\mleft(f\mright)\right),\]
and likewise for the fixed interpretation $r$ of the relation $\prec$,
\[\forall e, e', f, f' \in S_2 \left(e \equiv f \land e' \equiv f' \rightarrow \mleft\{r\mright\}\mleft(\tuple{e, e'}\mright) = \mleft\{r\mright\}\mleft(\tuple{f, f'}\mright)\right).\]
Thus the extensionality assertions $\left(7\right)$ and $\left(8\right)$ will be realised by any $\Sigma^1_1$ function with appropriate domain and range.

We have so far argued that assuming $\AC$, all the rules of inferences of intuitionistic logic, all logical axioms and all non-logical axioms of $\mathrm{CM}$ are realised. The following theorem thus follows from simple induction on the length of derivations:

\begin{theorem}[$\AC$]
    Whenever
    \[\mathrm{CM} \vdash \varphi\mleft(x_1, \ldots, x_n, X_1, \ldots, X_m, \mathbf{X}_1, \ldots, \mathbf{X}_r\mright),\]
    there exists $d \in \mathbb{N}$ such that $d \Vdash \varphi$.
\end{theorem}

\subsection{An ordinal analysis of \texorpdfstring{$\mathrm{CM}$}{CM}}

\newcommand{\ACi}{\text{$\Sigma^1_1$-$\mathrm{AC}^i$}}
\newcommand{\DC}{\text{$\Sigma^1_1$-$\mathrm{DC}$}}
\newcommand{\CA}{\text{$\Delta^1_1$-$\mathrm{CA}$}}

Let $\AC$ and $\DC$ denote the axioms in $\mathrm{ACA}$ plus the axioms of $\Sigma^1_1$-choice and $\Sigma^1_1$-dependent choice respectively, and let $\ACi$ denote the intuitionistic fragment of $\AC$. It is well-known that:

\begin{theorem}
    The following theories have the same proof-theoretic strength:
    \[\ACi \equiv \AC \equiv \DC.\]
\end{theorem}

Here, the proof-theoretic equivalence of $\ACi$ and $\AC$ is due to Aczel \cite{aczel77-martin-lof}. Friedman showed in \cite{friedman67-thesis} that one can interpret both $\AC$ and $\DC$ in $\CA$, and thus the three theories are equivalent in strength since $\CA$ is implied by $\AC$.

Since $\ACi$ is an obvious subtheory of $\mathrm{CM}$ and our analysis in previous sections shows that $\mathrm{CM}$ can be interpreted in $\AC$, it is an immediate corollary that:

\begin{corollary}
    $\mathrm{CM} \equiv \DC$. In particular, since the exact proof-theoretic strength of $\DC$ is analysed in \cite{friedman70-sigma11-12-ac}, we have $\left|\mathrm{CM}\right| = \left|\DC\right| = \varphi_{\varepsilon_0}\mleft(0\mright)$.
\end{corollary}

\section{Axiomatising a global well-ordering}
\label{sec:axiom-gwo}

\subsection{Axioms of a global well-ordering}

As mentioned in the introduction, in section 2.3 of \cite{weaver09-cm}, Weaver suggested briefly that one could extend the system $\mathrm{CM}$ by adding a ``genuine choice axiom'', stating that all the second-order objects in the theory are constructed in well-ordered stages. Formally, this can be represented as the following axioms:
\begin{enumerate}[label=$\left(\mathrm{W}\arabic*\right)$]
    \item the binary relation $\prec$ between second-order objects is a linear ordering;
    \item the transfinite induction scheme
          \[\forall X \left(\forall Y \left(Y \prec X \rightarrow \varphi\mleft(Y\mright)\right) \rightarrow \varphi\mleft(X\mright)\right) \rightarrow \forall X \, \varphi\mleft(X\mright),\]
          for a certain class of formulae;
    \item the axiom of countable initial segments
          \[\forall X \, \exists Z \, \forall Y \left(Y \prec X \rightarrow \exists n \, Y = \left(Z\right)_n\right).\]
\end{enumerate}

Do note that for the axiom scheme $\left(\mathrm{W}2\right)$ of transfinite induction, Weaver had been explicit that he did not accept the scheme for an arbitrary third-order formula $\varphi$ in the language, since it is obviously impredicative to assert that there exists a desired well-ordering $\prec$, on which transfinite induction holds for some formula $\varphi$ which involves the ordering relation $\prec$ itself. This is deemed unjustified according to his philosophical standpoint of mathematical conceptualism.

I plan to not delve deep into Weaver's philosophical arguments here, but one must point out that he failed to formally describe a class of formulae on which the axioms can characterise a restricted, predicative extension of $\mathrm{CM}$; and I believe it still remains open whether such a nice fragment of the full, unrestricted axiom scheme of transfinite induction can be found. More discussion on this will be reserved till the end in \autoref{sec:impredicativity}.

For now, we shall still stick to the entire collection of axioms $\left(\mathrm{W}1\right) + \left(\mathrm{W}2\right) + \left(\mathrm{W}3\right)$ where the formula $\varphi$ in the transfinite induction scheme ranges over all third-order formulae in the language, henceforth denoted $\mathrm{GWO}$. In \autoref{sec:cm-gwo-ordinal-analysis}, we will give an ordinal analysis of the extension $\mathrm{CM} + \mathrm{GWO}$, which will also serve as an upper bound of the proof-theoretic strength of any of its predicative fragments, if they actually exist.

\subsection{Zorn's lemma and applications}

Before the ordinal analysis, we shall note that while in \cite{weaver09-cm} Weaver gave a very detailed exposition of how ordinal mathematics, including theorems of real analysis, topology and functional analysis, can be treated inside the intuitionistic formal theory $\mathrm{CM}$, in order to justify its viability as a foundational theory of mathematics, he never did the same with this proposed extension, other than a short comment after proving the (separable) Hahn--Banach theorem:

\begin{theorem*}[$\mathrm{CM}$]
    Let $\mathbf{E}$ be a separable Banach space and $\mathbf{E}_0$ be a separable closed subspace. Then any bounded linear functional $\mathbf{f}_0 : \mathbf{E}_0 \rightarrow \mathbf{R}$ can be extended to a bounded linear functional $\mathbf{f} : \mathbf{E} \rightarrow \mathbf{R}$ satisfying $\left\|\mathbf{f}\right\| = \left\|\mathbf{f}_0\right\|$.
\end{theorem*}

Here, Weaver claimed that using the axioms for the global well-ordering $\prec$, this can be proved without assuming that $\mathbf{E}$ is separable through essentially the same proof.

In this section, we will expand on the direction by showing that a general result analogous to Zorn's lemma can be proved in $\mathrm{CM} + \mathrm{GWO}$; thus the global well-ordering axioms indeed serve their purpose as a choice principle and constitute indeed a useful extension to consider if one wants to use $\mathrm{CM}$ as a foundational theory.

To begin, we will first observe that omniscience holds for a bounded second-order quantifier:

\begin{lemma}[$\mathrm{CM} + \mathrm{GWO}$]
    \label{lem:bounded-second-order-omni}
    For any formula $\varphi\mleft(X\mright)$ and $\psi\mleft(X\mright)$, we have
    \[\forall X \left(\forall Y \prec X \left(\varphi\mleft(Y\mright) \lor \psi\mleft(Y\mright)\right) \rightarrow \left(\forall Y \prec X \, \varphi\mleft(Y\mright) \lor \exists Y \prec X \, \psi\mleft(Y\mright)\right)\right).\]
    Obviously the same also holds when every $\prec$ is replaced by $\preceq$ in the proposition.
\end{lemma}

\begin{proof}
    The axiom of countable initial segments asserts that for any second-order set $X$, there exists $Z$ satisfying
    \[\forall Y \prec X \, \exists n \, Y = \left(Z\right)_n.\]
    Assuming that $\forall Y \prec X \left(\varphi\mleft(Y\mright) \lor \psi\mleft(Y\mright)\right)$, then by numerical omniscience we have
    \[\forall n \left(\left(Z\right)_n \prec X \rightarrow \varphi\mleft(\left(Z\right)_n\mright)\right) \lor \exists n \left(\left(Z\right)_n \prec X \land \psi\mleft(\left(Z\right)_n\mright)\right).\]
    The former disjunct implies $\forall Y \prec X \, \varphi\mleft(Y\mright)$, while the latter disjunct implies $\exists Y \prec X \, \psi\mleft(Y\mright)$.
\end{proof}

We also need the fact that countable collections of second-order objects can be represented as one single second-order object in $\mathrm{CM}$. Instead of the obvious coding that uses slices $\left(Z\right)_n$, we will use a slightly modified construction so that we can accommodate the empty collection as well: for any second-order set $X$, let $\left[X\right]_n = \left\{m : \tuple{n, 0} \in_1 X \land \tuple{n, m + 1} \in_1 X\right\}$, and define a new binary relation $\in^*$ between \emph{two} second-order objects as
\[Y \in^* X \ \ \ratio\Rightarrow \ \ \exists n \left(\tuple{n, 0} \in_1 X \land Y = \left[X\right]_n\right).\]

We have some supplementary abbreviations such as
\begin{align*}
    X \subseteq^* Y \ \           & \ratio\Rightarrow \ \ \forall Z \in^* X \, Z \in^* Y,          \\
    X \subseteq^* \mathbf{X} \ \  & \ratio\Rightarrow \ \ \forall Z \in^* X \, Z \in_2 \mathbf{X}, \\
    X =^* Y \ \                   & \ratio\Rightarrow \ \ X \subseteq^* Y \land Y \subseteq^* X
\end{align*}
to compare these countable collections. One immediately observes that the omniscience statement for countable collections,
\[\forall Y \in^* X \left(\varphi\mleft(Y\mright) \lor \psi\mleft(Y\mright)\right) \rightarrow \left(\forall Y \in^* X \, \varphi\mleft(Y\mright) \lor \exists Y \in^* X \, \psi\mleft(Y\mright)\right),\]
follows directly from numerical omniscience, thus the relations defined above are all decidable.

In the same way, we can also use certain set-theoretic expressions over second-order objects as second-order terms. For example, we may have
\begin{align*}
    \left\{X_1, \ldots, X_k\right\}                    & \coloneqq \left\{\tuple{n, m} : 1 \leq n \leq k \land \left(m = 0 \lor \exists i \in_1 X_n \, m = i + 1\right)\right\},               \\
    \bigcup X                                          & \coloneqq \left\{\tuple{\tuple{n, m}, i} : \tuple{n, 0} \in_1 X \land \tuple{n, \tuple{m, i} + 1} \in_1 X\right\},                    \\
    \left\{Y \in^* X : \varphi\mleft(Y\mright)\right\} & \coloneqq \left\{\tuple{n, m} : \tuple{n, 0} \in_1 X \land \tuple{n, m} \in_1 X \land \varphi\mleft(\left[X\right]_n\mright)\right\},
\end{align*}
where in the last definition, the formula $\varphi\mleft(Y\mright)$ must be decidable. It is easy to verify that these collections ``contain'' the same second-order sets as their third-order counterparts.

We can now state our analogue of Zorn's lemma:

\begin{proposition}[$\mathrm{CM} + \mathrm{GWO}$]
    \label{prop:zorns-lemma}
    Fix a third-order set $\mathbf{A}$ and let $\varphi\mleft(X\mright)$ be a decidable formula (possibly with other parameters) on the countable subsets of $\mathbf{A}$, i.e.\ we have $\forall X \subseteq^* \mathbf{A} \left(\varphi\mleft(X\mright) \lor \neg \varphi\mleft(X\mright)\right)$. If $\varphi$ as a predicate is downward closed, i.e.
    \[\forall X, Y \left(\left(\varphi\mleft(X\mright) \land Y \subseteq^* X\right) \rightarrow \varphi\mleft(Y\mright)\right),\]
    and also closed under unions of countable chains, i.e.
    \[\forall X \left(\left(\forall Y, Z \in^* X \left(Y \subseteq^* Z \lor Z \subseteq^* Y\right) \land \forall Y \in^* X \, \varphi\mleft(Y\mright)\right) \rightarrow \varphi\mleft(\bigcup X\mright)\right),\]
    then there exists some $\mathbf{X} \subseteq \mathbf{A}$ satisfying $\forall X \subseteq^* \mathbf{X} \ \varphi\mleft(X\mright)$, and is maximal among such subsets in the sense that
    \[\forall Y \left(Y \not\in \mathbf{X} \rightarrow \exists Z \subseteq^* \mathbf{X} \ \neg \varphi\mleft(Z \cup \left\{Y\right\}\mright)\right).\]
\end{proposition}

\begin{proof}
    Let
    \begin{align*}
        \psi\mleft(X, Z\mright) \ \    & \coloneqq \ \ \varphi\mleft(Z\mright) \land \left(\varphi\mleft(Z \cup \left\{X\right\}\mright) \rightarrow X \in^* Z\right),             \\
        \theta\mleft(X, Z\mright) \ \  & \coloneqq \ \ \forall Y \in^* Z \, Y \preceq X \land \forall Y \preceq X \, \psi\mleft(Y, \left\{W \in^* Z : W \preceq Y\right\}\mright).
    \end{align*}
    We first show by transfinite induction that
    \[\forall X, Z_1, Z_2 \left(\left(\theta\mleft(X, Z_1\mright) \land \theta\mleft(X, Z_2\mright)\right) \rightarrow Z_1 =^* Z_2\right).\]
    This is because for any $Y \prec X$, $\theta\mleft(X, Z_1\mright)$ implies $\theta\mleft(Y, \left\{W \in^* Z_1 : W \preceq Y\right\}\mright)$ and likewise for $Z_2$. Therefore, from the inductive hypothesis we know that
    \[\left\{W \in^* Z_1 : W \preceq Y\right\} =^* \left\{W \in^* Z_2 : W \preceq Y\right\},\]
    which implies $Y \in^* Z_1 \leftrightarrow Y \in^* Z_2$. Now assume that $X \in^* Z_1$, then $Z_1 =^* Z_2 \cup \left\{X\right\}$, so $\varphi\mleft(Z_2 \cup \left\{X\right\}\mright)$ holds and thus $X \in^* Z_2$. The same can be said assuming $X \in^* Z_2$, so we can finally conclude that $Z_1 =^* Z_2$.

    We then show by transfinite induction again that $\forall X \, \exists Z \, \theta\mleft(X, Z\mright)$. For any $X$, let $W$ be such that $\forall Y \prec X \, \exists n \, Y = \left(W\right)_n$, and the inductive hypothesis implies
    \[\forall n \, \exists Z \left(\left(W\right)_n \prec X \rightarrow \theta\mleft(\left(W\right)_n, Z\mright)\right).\]
    By dependent choice, we find some $U$ satisfying $\forall n \left(\left(W\right)_n \prec X \rightarrow \theta\mleft(\left(W\right)_n, \left(U\right)_n\mright)\right)$. Now define
    \[Z' = \left\{\tuple{n, m} : \left(W\right)_n \prec X \land \left(m = 0 \lor \exists i \in_1 \left(U\right)_n \, m = i + 1\right)\right\},\]
    then for each $Y \in^* Z'$, we have $\varphi\mleft(Y\mright)$. Here, $Z'$ is additionally a chain because for any $\left[Z'\right]_m, \left[Z'\right]_n \in^* Z'$, we can assume without loss of generality that $\left(W\right)_m \preceq \left(W\right)_n$, and we will have
    \[\left[Z'\right]_m =^* \left\{W \in^* \left[Z'\right]_n : W \preceq \left(W\right)_m\right\} \subseteq^* \left[Z'\right]_n\]
    from the previous induction argument. Hence $\varphi\mleft(\bigcup Z'\mright)$. It also follows that $\left[Z'\right]_n = \left\{Y \in^* \bigcup Z' : Y \prec \left(W\right)_n\right\}$ for any $n$ such that $\left(W\right)_n \prec X$.

    Therefore, it suffices to set $Z = \left(\bigcup Z'\right) \cup \left\{X\right\}$ if $\varphi\mleft(\left(\bigcup Z'\right) \cup \left\{X\right\}\mright)$ holds and $Z = \bigcup Z'$ otherwise (which is possible because $\varphi$ is decidable). It is easy to verify that $\theta\mleft(X, Z\mright)$ holds indeed.

    Finally, we observe that $\theta\mleft(X, Z\mright)$ is decidable by \autoref{lem:bounded-second-order-omni}. Therefore, we can comprehend the diagonal set
    \[\mathbf{X} = \left\{Y : \exists Z \left(\theta\mleft(Y, Z\mright) \land Y \in^* Z\right)\right\}.\]
    For any countable subset $X \subseteq^* \mathbf{X}$, we have
    \[\forall n \, \exists W \left(\tuple{n, 0} \in_1 X \rightarrow \forall Y \preceq \left[X\right]_n \, \exists m \, Y = \left(W\right)_m\right),\]
    and by dependent choice, we can find some $U$ satisfying
    \[\forall n \left(\tuple{n, 0} \in_1 X \rightarrow \forall Y \preceq \left[X\right]_n \, \exists m \, Y = \left(\left(U\right)_n\right)_m\right).\]
    We can now diagonalise and construct
    \[Y_0 = \left\{\tuple{n, m} : \tuple{n, m} \not\in_1 \left(\left(U\right)_n\right)_m\right\}\]
    so that $\forall Y \in^* X \, \forall Y' \preceq Y \ Y' \prec Y_0$. Let $Z_0$ be a countable collection satisfying $\theta\mleft(Y_0, Z_0\mright)$, then our definition of $\mathbf{X}$ and the uniqueness of $Z_0$ as shown previously ensure that for any $Z \prec Y_0$,
    \[Y \in_2 \mathbf{X} \leftrightarrow Y \in^* Z_0.\]
    In other words, $X \subseteq^* Z_0$, and $\varphi\mleft(X\mright)$ follows immediately from $\varphi\mleft(Z_0\mright)$.

    On the other hand, suppose that some $Y \not\in_2 \mathbf{X}$. Likewise, there exists $Z$ satisfying $\theta\mleft(Y, Z\mright)$ and $Y \not\in^* Z$, which, by our stipulation of $\psi$ at the beginning of this proof, imply $\neg \varphi\mleft(Z \cup \left\{Y\right\}\mright)$. The uniqueness of $Z$ ensures that also $Z \subseteq^* \mathbf{X}$. Therefore, we have checked that $\mathbf{X}$ has all the desired properties.
\end{proof}

This can be applied just like the classical Zorn's lemma, to show that one can find maximal subsets satisfying certain properties in uncountable collections, as long as the corresponding decidability criterion can be verified through intuitionistic reasoning. As an example, let us consider third-order objects $\mathbf{F}$, $\mathbf{V}$ such that $\mathbf{F}$ is a field and $\mathbf{V}$ is an $\mathbf{F}$-vector space, with their respective unary and binary operations also given as third-order objects, following the conventions in section 3 of \cite{weaver09-cm}. We can derive that

\begin{corollary}[$\mathrm{CM} + \mathrm{GWO}$]
    Suppose that the formula
    \[\varphi\mleft(k, X\mright) \quad \coloneqq \quad \exists Y \left(\forall n \leq k \left(Y\right)_n \in \mathbf{F} \land \exists n \leq k \left(Y\right)_n \neq 0_\mathbf{F} \land 0_\mathbf{V} = \sum_{n = 0}^k \left(Y\right)_n \cdot \left(X\right)_n\right),\]
    i.e.\ informally the statement ``$\left(X\right)_0, \ldots, \left(X\right)_k$ are linearly dependent as vectors in $\mathbf{V}$'', is decidable (whenever $\left(X\right)_0, \ldots, \left(X\right)_k \in_2 \mathbf{V}$), then $\mathbf{V}$ has a basis $\mathbf{X} \subseteq \mathbf{V}$ such that every vector $V \in \mathbf{V}$ is spanned by finitely many elements in $\mathbf{X}$.
\end{corollary}

\begin{proof}
    Fix a way to code finite tuples of natural numbers, so that every $a \in \mathbb{N}$ codes some tuple $\tuple{\left(a\right)_0, \ldots, \left(a\right)_{\left|a\right|}}$. Now, let
    \[\psi\mleft(X\mright) \quad \coloneqq \quad X \subseteq^* \mathbf{V} \land \forall a \left(\exists n \leq \left|a\right| \tuple{\left(a\right)_n, 0} \not\in_1 X \lor \neg \varphi\mleft(\left|a\right|, T\mleft(X, a\mright)\mright)\right)\]
    where the tuple $T\mleft(X, a\mright) = \left\{\tuple{n, m} : n \leq \left|a\right| \land \tuple{\left(a\right)_n, m + 1} \in X\right\}$ collects precisely the elements $\left[X\right]_{\left(a\right)_0}, \ldots, \left[X\right]_{\left(a\right)_{\left|a\right|}}$. By numerical omniscience, $\psi$ is also decidable.

    It is easy to check that the premises of \autoref{prop:zorns-lemma} are all satisfied, thus there is a maximal set $\mathbf{X}$ satisfying $\forall X \subseteq^* \mathbf{X} \, \psi\mleft(X\mright)$, which also implies $\mathbf{X} \subseteq \mathbf{V}$. Additionally, for any vector $V \in_2 \mathbf{V}$, if $V \in_2 \mathbf{X}$, then $V$ is trivially spanned; otherwise $V \not\in_2 \mathbf{X}$, and we know from decidability that $\neg \psi\mleft(Z \cup \left\{V\right\}\mright)$ for some $Z \subseteq^* \mathbf{X}$, i.e.\ a finite tuple of vectors from $Z \cup \left\{V\right\}$ are linearly dependent. Since $Z$ alone does not span $0_\mathbf{V}$, $V$ must be spanned by elements in $Z$.
\end{proof}

In general, the criterion $\varphi\mleft(k, X\mright)$ may not be decidable for every vector space over every field, due to the second-order quantifier in the beginning. However, very often the field or the vector space in concern is nice enough for the basis to exist:
\begin{itemize}
    \item When \emph{the field $\mathbf{F}$ is countable}, $\varphi\mleft(k, X\mright)$ is decidable simply by numerical omniscience. This shows that, for example, the Hamel basis for the reals $\mathbf{R}$ as a vector space over the rationals $\mathbf{Q}$ exists.
    \item When \emph{the vector space $\mathbf{V}$ has an inner product $\tuple{\cdot, \cdot} : \mathbf{V} \times \mathbf{V} \rightarrow \mathbf{F}$}, given some finite collection $X$ of vectors in $\mathbf{V}$, we can apply the Gram--Schmidt process by computing
          \begin{align*}
              \left(Y\right)_n & = \left(X\right)_n - \sum_{m = 0}^{n - 1} \tuple{\left(Z\right)_m, \left(X\right)_n} \left(Z\right)_m, \\
              \left(Z\right)_n & = \tuple{\left(Y\right)_n, \left(Y\right)_n}^{-1} \left(Y\right)_n,
          \end{align*}
          for $0 \leq n \leq k$ and halt when we encounter some $\left(Y\right)_n = 0_\mathbf{V}$. Obviously $\exists n \leq k \left(Y\right)_n = 0_\mathbf{V}$ implies $\varphi\mleft(k, X\mright)$; and conversely if we have $0_\mathbf{V} = \sum_{n = 0}^k \left(W\right)_n \cdot \left(X\right)_n$ where $\left(W\right)_n \in_2 \mathbf{F}$ are not all zero, we similarly also have the equality
          \[\left(X\right)_n = \sum_{m = 0}^{n - 1} - \left(W\right)_m \left(W\right)_n^{-1} \left(X\right)_m\]
          for some $0 \leq n \leq k$ and thus $\left(Y\right)_n = 0_\mathbf{V}$ assuming $\left(Y\right)_m \neq 0_\mathbf{V}$ for all $m < n$. In this case, $\varphi\mleft(k, X\mright)$ is equivalent to $\exists n \leq k \left(Y\right)_n = 0_\mathbf{V}$ and thus decidable.
\end{itemize}

We shall sketch that $\varphi\mleft(k, X\mright)$ is decidable also when \emph{$\mathbf{F}$ is a complete $\mathbf{R}$-normed field with a countable dense set $C \subseteq^* \mathbf{F}$ with respect to the induced norm}, and \emph{$\mathbf{V}$ has a compatible norm $\left\|\cdot\right\| : \mathbf{V} \rightarrow \mathbf{R}$}. Any real or complex Banach space, for example, falls in this category.

\begin{corollary}
    Let $\mathbf{F}$ be a complete $\mathbf{R}$-normed field with a countable dense set $C \subseteq^* \mathbf{F}$ with respect to the induced norm, and $\mathbf{V}$ be a $\mathbf{F}$-vector space with a compatible norm $\left\|\cdot\right\| : \mathbf{V} \rightarrow \mathbf{R}$. Then $\mathbf{V}$ has a basis.
\end{corollary}

\begin{proof}
    We consider the condition that
    \[\psi\mleft(k, X\mright) \ \ \coloneqq \ \ \exists n \leq k \, \forall z \, \exists m^{\left(z\right)}_0, \ldots, m^{\left(z\right)}_{n - 1} \left\|\left(X\right)_n - \sum_{i = 0}^{n - 1} \left[C\right]_{m^{\left(z\right)}_i} \left(X\right)_i\right\| < 1 / z,\]
    which is decidable by numerical omniscience.

    First observe that if we have some equality $0_\mathbf{V} = \sum_{n = 0}^k \left(Y\right)_n \cdot \left(X\right)_n$, by numerical omniscience we can find the last $n \leq k$ such that $\left(Y\right)_n \neq 0_\mathbf{F}$ and rewrite it as
    \[\left(X\right)_n = \sum_{m = 0}^{n - 1} - \left(Y\right)_m \left(Y\right)_n^{-1} \left(X\right)_m.\]
    The density of $C$ ensures that we can find elements in $C$ close enough to each scalar $- \left(Y\right)_m \left(Y\right)_n^{-1} \in_2 \mathbf{F}$ so that $\psi\mleft(k, X\mright)$ holds.

    Now assume instead that $\psi\mleft(k, X\mright)$ holds, then by its decidability and induction on the natural numbers we can find the least such $n \leq k$. By dependent choice, fix any $0 \leq i < n$ and we can collect the sequence of $\left[C\right]_{m^{\left(z\right)}_i}$ used in $\psi$ for every $z$. It suffices to show that each such sequence is Cauchy, because then the limits would describe $\left(X\right)_n$ as a linear combination of $\left(X\right)_0, \ldots, \left(X\right)_{n - 1}$.

    Here, the minimality of $n$ will ensure that $\left(X\right)_0, \ldots, \left(X\right)_{n - 1}$ are linearly independent, i.e.\ the function
    \[\mathbf{f} : \tuple{Y_0, \ldots, Y_{n - 1}} \in \mathbf{F}^n \mapsto \sum_{m = 0}^{n - 1} Y_m \cdot \left(X\right)_m\]
    will be injective. With respect to the $1$-norm
    \[\left\|\tuple{Y_0, \ldots, Y_{n - 1}}\right\| = \sum_{m = 0}^{n - 1} \left\|Y_m\right\|\]
    on $\mathbf{F}^n$, $\mathbf{f}$ is continuous. Thus, if we consider the separable compact subspace $\mathbf{S}^1 = \left\{Y \in_2 \mathbf{F}^n : \left\|Y\right\| = 1\right\}$, then by proposition 3.40 in \cite{weaver09-cm}, the image $\mathbf{f}\mleft(\mathbf{S}^1\mright)$ exists and is compact in $\mathbf{V}$. Since $\left(X\right)_0, \ldots, \left(X\right)_{n - 1}$ are linearly independent, $0_\mathbf{V} \not\in_2 \mathbf{f}\mleft(\mathbf{S}^1\mright)$, thus there is a minimum $M > 0$ such that $\forall V \in_2 \mathbf{f}\mleft(\mathbf{S}^1\mright) \left\|V\right\| \geq M$. Now, for any $Y \in_2 \mathbf{F}^n$, $Y \neq 0_{\mathbf{F}^n}$ such that $\left\|\mathbf{f}\mleft(Y\mright)\right\| < \varepsilon$, we have
    \[M \leq \left\|\mathbf{f}\mleft(\left\|Y\right\|^{-1} Y\mright)\right\| < \left\|Y\right\|^{-1} \varepsilon,\]
    i.e.\ each $\left\|Y_i\right\| < M^{-1} \varepsilon$. This justifies precisely what we need, that when the sum $\sum_{i = 0}^{n - 1} \left[C\right]_{m^{\left(z\right)}_i} \left(X\right)_i$ is Cauchy in the parameter $z$, then so is the sequence $\left[C\right]_{m^{\left(z\right)}_i}$ for every fixed $i$.
\end{proof}

\section{An ordinal analysis of \texorpdfstring{$\mathrm{CM} + \mathrm{GWO}$}{CM + GWO}}
\label{sec:cm-gwo-ordinal-analysis}

\subsection{Hyperarithmetic theory in \texorpdfstring{$\mathrm{ATR}$}{ATR}}
\label{subsec:hyperarithmetic}

Now, we will return from our detour on Zorn's lemma back to providing an ordinal analysis of the extended theory $\mathrm{CM} + \mathrm{GWO}$. We shall do so by reusing the realisability model defined in \autoref{subsec:realisability-cm}, but move to a stronger base theory that provides a nice well-ordering on the collection of second-order realisers $S_2$. This well-ordering will be generated from the hierarchy of hyperarithmetic sets given by iterating Turing jumps.

For now, let us work in the theory $\mathrm{ATR}$ (with unrestricted induction). Everything we have done so far using $\Sigma^1_1$ functions apply immediately due to the following result, which is theorem V.8.3 in Simpson \cite{simpson09-soa}:

\begin{fact*}
    $\mathrm{ATR}$ implies $\AC$.
\end{fact*}

The constructions we will introduce here roughly follow section VIII.3 in \cite{simpson09-soa}: we again start from a universal lightface $\Pi^0_1$-formula $\pi\mleft(e, m, X\mright)$ with all free variables indicated, that is, for any $\Pi^0_1$-formula $\varphi\mleft(e', m, X\mright)$ with the same number of free variables, $\mathrm{ACA}_0$ proves that
\[\forall e' \, \exists e \, \forall m \, \forall X \left(\varphi\mleft(e', m, X\mright) \leftrightarrow \pi\mleft(e, m, X\mright)\right).\]
We shall make use of the following series of definitions:

\begin{definition}[Relative recursiveness]
    We say that a set $Y \subseteq \mathbb{N}$ is \emph{$X$-recursive} and write $Y \leq_T X$ if there exists $e = \tuple{e_0, e_1} \in \mathbb{N}$ such that
    \[Y = \left\{m : \pi\mleft(e_0, m, X\mright)\right\} = \left\{m : \neg \pi\mleft(e_1, m, X\mright)\right\}.\]
    We say that $e$ is the \emph{$X$-recursive index} of $Y$. We write
    \[\Rec^X = \left\{e \in \mathbb{N} : \forall m \left(\pi\mleft(e_0, m, X\mright) \leftrightarrow \neg \pi\mleft(e_1, m, X\mright)\right)\right\}\]
    to denote the set of all $X$-recursive indices, and for any $e \in \Rec^X$, let $R^X_e$ denote the set defined by $e$. The superscript $X$ here will be omitted when $X = \varnothing$.
\end{definition}

\begin{definition}[Ordinal notations]
    Let $\mathcal{O}_+$ denote the set of indices $\alpha \in \mathbb{N}$ where $\alpha_0 \in \Rec$, $R_{\alpha_0}$ is a linear ordering on a set $\field\mleft(R_{\alpha_0}\mright) \subseteq \mathbb{N}$ and $\alpha_1 \in \field\mleft(R_{\alpha_0}\mright)$. The canonical partial ordering on $\mathcal{O}_+$ is given by
    \[\alpha <_\mathcal{O} \beta \ \ \ratio\Leftrightarrow \ \ \alpha_0 = \beta_0 \land \tuple{\alpha_1, \beta_1} \in R_{\alpha_0}.\]

    Let
    \[\mathcal{O} = \left\{\alpha \in \mathcal{O}_+ : \neg \exists f : \mathbb{N} \rightarrow \mathcal{O}_+ \left(f\mleft(0\mright) = \alpha \land \forall n \, f\mleft(n + 1\mright) <_\mathcal{O} f\mleft(n\mright)\right)\right\},\]
    i.e.\ denote the subclass of indices $\alpha$ such that $<_\mathcal{O}$ is a well-ordering on the initial segment before $\alpha$. Given any $\alpha \in \mathcal{O}$, we say its \emph{order type}\footnote{In our setup, we will not define concrete representations for order types of well-orderings. Instead, we will treat statements $\otyp\mleft(X\mright) = \otyp\mleft(Y\mright)$ and $\otyp\mleft(X\mright) \leq \otyp\mleft(Y\mright)$ as abbreviations for the assertions that corresponding order-preserving injections/bijectons exist respectively.} is precisely the order type of this initial segment
    \[\otyp\mleft(\alpha\mright) = \otyp\mleft(\left\{\beta \in \mathcal{O}_+ : \beta <_\mathcal{O} \alpha\right\}\mright).\]
\end{definition}

\begin{definition}[$\mathrm{H}$-sets]
    The \emph{Turing jump} of a set $X \subseteq \mathbb{N}$ is given by
    \[\TJ\mleft(X\mright) = \left\{\tuple{e, m} : \pi\mleft(e, m, X\mright)\right\}.\]
    The following formula identifies when $Y$ represents the iteration of Turing jumps along some ordering $R_{\alpha_0}$ below $\alpha_1$:
    \[H\mleft(\alpha, Y\mright) \ \ \ratio\Leftrightarrow \ \ \alpha \in \mathcal{O}_+ \land Y = \left(Y\right)_{<_\mathcal{O} \alpha} \land \forall \beta <_\mathcal{O} \alpha \, \left(Y\right)_\beta = \TJ\mleft(\left(Y\right)_{<_\mathcal{O} \beta}\mright),\]
    where the truncation of $Y$ before an ordinal notation $\alpha$ is
    \[\left(Y\right)_{<_\mathcal{O} \alpha} = \left\{\tuple{\beta, m} : \beta <_\mathcal{O} \alpha \land m \in \left(Y\right)_{\beta}\right\}.\]
\end{definition}

By theorem VIII.3.15 in \cite{simpson09-soa}, $\mathrm{ATR}_0$ proves $\forall \alpha \in \mathcal{O} \ \exists Y \, H\mleft(\alpha, Y\mright)$. We can thus make use of the following property, which is lemma VIII.3.14 in \cite{simpson09-soa}:

\begin{proposition}[$\mathrm{ACA}_0$]
    Given $\alpha, \beta \in \mathcal{O}$ and some $Y$ such that $H\mleft(\alpha, Y\mright)$, then the order types $\otyp\mleft(\alpha\mright)$ and $\otyp\mleft(\beta\mright)$ are comparable, and the comparison map is $\leq_T Y$.
\end{proposition}

\begin{corollary}[$\mathrm{ATR}_0$]
    \label{cor:well-order-cmp}
    There exists $\mathrm{cmp} \in \left[\mathcal{O}, S_2\right]$ such that
    \[\mleft\{\mleft\{\mathrm{cmp}\mright\}\mleft(\alpha\mright)\mright\}\mleft(\beta\mright) = 1 \ \ \Leftrightarrow \ \ \beta \in \mathcal{O} \land \otyp\mleft(\beta\mright) \leq \otyp\mleft(\alpha\mright).\]
\end{corollary}

\begin{proof}
    Given $\alpha \in \mathcal{O}$, fix $Y$ to be a set such that $H\mleft(\alpha, Y\mright)$. We let $\mleft\{\mleft\{\mathrm{cmp}\mright\}\mleft(\alpha\mright)\mright\}\mleft(\beta\mright)$ be $1$ if $\beta \in \mathcal{O}_+$ and there exists $e \in \Rec^Y$ such that $R^Y_e$ is an order-isomorphism between $\left\{\gamma \in \mathcal{O}_+ : \gamma <_\mathcal{O} \beta\right\}$ and an initial segment of $\left\{\gamma \in \mathcal{O}_+ : \gamma <_\mathcal{O} \alpha\right\}$. Let $\mleft\{\mleft\{\mathrm{cmp}\mright\}\mleft(\alpha\mright)\mright\}\mleft(\beta\mright)$ be $0$ otherwise. Notice that this statement is $\Sigma^1_1$.

    Obviously, suppose that $\mleft\{\mleft\{\mathrm{cmp}\mright\}\mleft(\alpha\mright)\mright\}\mleft(\beta\mright) = 1$, then $\beta \in \mathcal{O}$ because any strictly decreasing sequence below $\beta$ will correspond to a strictly decreasing sequence below $\alpha$ through $R^Y_e$. Conversely, suppose that indeed $\beta \in \mathcal{O}$ and $\otyp\mleft(\beta\mright) \leq \otyp\mleft(\alpha\mright)$, then the proposition above ensures that the comparison function is $Y$-recursive, so $\mleft\{\mleft\{\mathrm{cmp}\mright\}\mleft(\alpha\mright)\mright\}\mleft(\beta\mright)$ must be $1$.
\end{proof}

\subsection{Code conversion}

Denote
\[\mathcal{H} = \left\{a : a_0 \in \mathcal{O} \land \exists Y \left(H\mleft(a_0, Y\mright) \land a_1 \in \Rec^Y\right)\right\}.\]
We say that a class $X \subseteq \mathbb{N}$ is \emph{hyperarithmetic}, denoted $X \leq_H \varnothing$, if there exists $a \in \mathcal{H}$ and $Y \subseteq \mathbb{N}$ such that $H\mleft(a_0, Y\mright)$ and $X = R^Y_{a_1}$. We say that $a = \tuple{a_0, a_1}$ is a \emph{hyperarithmetic index} of $X$.

The aim of this section is to show that one can convert between the hyperarithmetic index of a set and the $\Sigma^1_1$ function code of its characteristic function. Namely, we will prove the following main theorem:

\begin{theorem}[$\mathrm{ATR}_0$]
    \label{thm:conv-s2-to-hyp}
    There exists $\mathrm{hyp} \in \left[S_2, \mathcal{H}\right]$ such that
    \begin{itemize}
        \item for each $e \in S_2$, $\mleft\{\mathrm{hyp}\mright\}\mleft(e\mright)$ is a hyperarithmetic index of $\left\{n : \mleft\{e\mright\}\mleft(n\mright) = 1\right\}$;
        \item for $d, e \in S_2$, $\mleft\{\mathrm{hyp}\mright\}\mleft(d\mright) = \mleft\{\mathrm{hyp}\mright\}\mleft(e\mright)$ if and only if $d \equiv e$.
    \end{itemize}
\end{theorem}

\begin{proof}
    To begin, we will replicate some constructions in the proof of theorem VIII.3.19 in \cite{simpson09-soa} on the set $\left\{n : \mleft\{e\mright\}\mleft(n\mright) = 1\right\}$: observe that the formula
    \[\mleft\{e\mright\}\mleft(n\mright) = 1 \ \ \Leftrightarrow \ \ \mleft\{e\mright\}\mleft(n\mright) \neq 0 \ \ \Leftrightarrow \ \ \forall X \, \neg \pi\mleft(e_0, e_1, n, 0, X\mright)\]
    is $\Pi^1_1$, for any fixed parameter $e \in S_2$. By Kleene normal form theorem (theorem V.1.4 in \cite{simpson09-soa}), there exists a $\Sigma^0_0$-formula $\theta\mleft(e, n, \tau\mright)$ such that $\mathrm{ACA}_0$ proves
    \[\forall n \left(\mleft\{e\mright\}\mleft(n\mright) = 1 \leftrightarrow \forall f \, \exists k \, \neg \theta\mleft(e, n, f\mleft[k\mright]\mright)\right)\]
    where $f$ ranges over second-order objects representing functions $\mathbb{N} \rightarrow \mathbb{N}$, and $f\mleft[k\mright] = \tuple{f\mleft(0\mright), \ldots, f\mleft(k - 1\mright)} \in \mathbb{N}$. We then consider recursive trees
    \[T^{(e)}_n = \left\{\tau \in \Seq : \forall k \leq \mathrm{len}\mleft(\tau\mright) \ \theta\mleft(e, n, \tau\mleft[k\mright]\mright)\right\}\]
    and let $\KB^{(e)}_n = \KB\mleft(T^{(e)}_n\mright)$ be the Kleene-Brouwer ordering of $T^{(e)}_n$. Following the discussions in section V.1 in \cite{simpson09-soa}, we have
    \begin{align*}
        \mleft\{e\mright\}\mleft(n\mright) = 1 \ \  & \Leftrightarrow \ \ \forall f \, \exists k \, \neg \theta\mleft(e, n, f\mleft[k\mright]\mright) \\
                                                    & \Leftrightarrow \ \ T^{(e)}_n \text{ has no path}                                               \\
                                                    & \Leftrightarrow \ \ \KB^{(e)}_n \text{ is a recursive well-ordering}.
    \end{align*}

    Now, as mentioned in \autoref{subsec:hyperarithmetic}, we can assume that $\forall \alpha \in \mathcal{O} \ \exists Y \, H\mleft(\alpha, Y\mright)$. Then Simpson's argument applies to justify that the collection of well-orderings $\KB^{(e)}_n$ (for all $n \in \mathbb{N}$, given a fixed $e$) is bounded in $\mathcal{O}$. Let $\alpha \in \mathcal{O}$ be such a bound and $\Omega \subseteq \mathbb{N}$ be such that $H\mleft(\alpha, \Omega\mright)$. By lemma VIII.3.14 in \cite{simpson09-soa}, the order types of any initial segment of $\KB^{(e)}_n$ and any $\beta \leq_\mathcal{O} \alpha$ will all be comparable via comparison maps $\leq_T \Omega$, which means for any such well-orderings $A, B$, a comparison statement $\otyp\mleft(A\mright) \leq \otyp\mleft(B\mright)$ will be an arithmetic formulae ``there exists $a \in \mathbb{N}$ such that $R^\Omega_a$ is an order isomorphism between $A$ and an initial segment of $B$''.

    Therefore, by $\Delta^1_1$-comprehension,
    \[\Gamma = \left\{\beta \leq_\mathcal{O} \alpha : \forall n \left(\mleft\{e\mright\}\mleft(n\mright) = 1 \rightarrow \otyp\mleft(\KB^{(e)}_n\mright) \leq \otyp\mleft(\beta\mright)\right)\right\}\]
    is a set and we can find $\lambda = \min \Gamma$. It is easy to also observe that for any $\gamma <_\mathcal{O} \lambda$, there must exist $n \in \mathbb{N}$ such that $\mleft\{e\mright\}\mleft(n\mright) = 1$ and $\otyp\mleft(\gamma\mright) < \otyp\mleft(\KB^{(e)}_n\mright)$. In fact, these two characteristics determine $\lambda$ uniquely up to order isomorphism, i.e.\ we have the lemma that

    \begin{lemma*}
        Suppose that $\rho \in \mathcal{O}_+$ satisfies both $\otyp\mleft(\KB^{(e)}_n\mright) \leq \otyp\mleft(\rho\mright)$ for any $n$ such that $\mleft\{e\mright\}\mleft(n\mright) = 1$, and also for any $\gamma <_\mathcal{O} \rho$, there must exist $n \in \mathbb{N}$ such that $\mleft\{e\mright\}\mleft(n\mright) = 1$ and $\otyp\mleft(\gamma\mright) < \otyp\mleft(\KB^{(e)}_n\mright)$, then $\otyp\mleft(\rho\mright) = \otyp\mleft(\lambda\mright)$.
    \end{lemma*}

    \begin{proof}
        Firstly, assume that $\rho \not\in \mathcal{O}$, then there exists a strictly decreasing sequence $f : \mathbb{N} \rightarrow \mathcal{O}_+$ such that $f\mleft(0\mright) = \rho$. We have $f\mleft(1\mright) <_\mathcal{O} \rho$, hence we can consider some order-preserving map $g : \left\{\gamma : \gamma <_\mathcal{O} f\mleft(1\mright)\right\} \rightarrow \KB^{(e)}_n$. Now $m \mapsto g\mleft(f\mleft(m + 2\mright)\mright)$ is a descending sequence in $\KB^{(e)}_n$, a contradiction.

        Therefore, we must have $\rho \in \mathcal{O}$. We then know that the order types of $\rho$ and $\lambda$ are comparable, but it is easy to verify that we can have neither $\otyp\mleft(\rho\mright) > \otyp\mleft(\lambda\mright)$ nor $\otyp\mleft(\rho\mright) < \otyp\mleft(\lambda\mright)$.
    \end{proof}

    Finally, given that $\lambda \in \mathcal{O}$, we can remove the ambiguity around order-isomorphic codes by requiring that $\lambda$ is minimal among all order-isomorphic elements in $\mathcal{O}$, with respect to the natural ordering of $\mathbb{N}$. Such a minimum exists by $\Sigma^0_0$-induction because, using some $Y \subseteq \mathbb{N}$ satisfying $H\mleft(\lambda, Y\mright)$ as a parameter, the set of elements in $\mathcal{O}$ order-isomorphic to $\lambda$ is an arithmetic set. The rest of Simpson's argument in VIII.3.19 applies to show that if $H\mleft(\lambda + \omega, Y\mright)$, then $\left\{n : \mleft\{e\mright\}\mleft(n\mright) = 1\right\} \leq_T Y$.

    Therefore, we can stipulate $\mleft\{\mathrm{hyp}\mright\}\mleft(e\mright)$ to compute some $a \in \mathbb{N}$ such that there exists $\alpha \in \mathcal{O}_+$ and $Y \subseteq \mathbb{N}$ and we have
    \begin{itemize}
        \item $H\mleft(a_0, Y\mright)$ and\footnote{Here we fix $R_{\alpha + \omega}$ to mean the following specific recursive well-ordering: simply order the set $\left\{2\beta + 2 : \beta <_\mathcal{O} \alpha\right\}$ by
                  \[2\beta_1 + 2 < 2\beta_2 + 2 \ \ \ratio\Leftrightarrow \ \ \beta_1 <_\mathcal{O} \beta_2,\]
                  and append the sequence $1, 3, 5, 7, \ldots, 0$. Then the initial segment before $0$ has order type $\alpha + \omega$.} $R_{a_0} = R_{\alpha + \omega}$;
        \item for any $n \in \mathbb{N}$, either $\mleft\{e\mright\}\mleft(n\mright) = 0$ or $\otyp\mleft(\KB^{(e)}_n\mright) \leq \otyp\mleft(\alpha\mright)$;
        \item for any $\beta <_\mathcal{O} \alpha$, there exists $n \in \mathbb{N}$ such that $\mleft\{e\mright\}\mleft(n\mright) = 1$ and $\otyp\mleft(\beta\mright) < \otyp\mleft(\KB^{(e)}_n\mright)$;
        \item for any $\gamma \in \mathcal{O}_+$, if $\gamma$ is order-isomorphic to $a_0$, then $\gamma \geq a_0$;
        \item $a_1 \in \Rec^Y$ and for any $n \in \mathbb{N}$, $\mleft\{e\mright\}\mleft(n\mright) = 1 \lor n \not\in R^Y_{a_1}$ and $\mleft\{e\mright\}\mleft(n\mright) = 0 \lor n \in R^Y_{a_1}$;
        \item for any $b \in \Rec^Y$, if $n \in R^Y_{a_1} \leftrightarrow n \in R^Y_b$ holds for any $n \in \mathbb{N}$, then $b \geq a_1$.
    \end{itemize}
    This determines $a$ uniquely and is $\Sigma^1_1$. Since our analysis above shows that such $a$ must exist for any $e \in S_2$, we have $S_2 \subseteq \dom\mleft(\mathrm{hyp}\mright)$.

    To verify that $\mleft\{\mathrm{hyp}\mright\}\mleft(d\mright) = \mleft\{\mathrm{hyp}\mright\}\mleft(e\mright)$ if and only if $d \equiv e$, simply observe that $\mleft\{\mathrm{hyp}\mright\}\mleft(e\mright)$ is uniquely determined by the set $\left\{n : \mleft\{e\mright\}\mleft(n\mright) = 1\right\}$.
\end{proof}

A (rough) inverse of the conversion is very easy to compute:

\begin{proposition}[$\mathrm{ACA}_0$]
    There exists $\mathrm{inv} \in \left[\mathcal{H}, S_2\right]$ such that, for any $a \in \mathcal{H}$ and any $Y \subseteq \mathbb{N}$ satisfying $H\mleft(a_0, Y\mright)$,
    \[\mleft\{\mleft\{\mathrm{inv}\mright\}\mleft(a\mright)\mright\}\mleft(n\mright) = 1 \ \ \Leftrightarrow \ \ n \in R^Y_{a_1}.\]
\end{proposition}

\begin{proof}
    This is simply because the formula
    \[\mleft\{\mleft\{\mathrm{inv}\mright\}\mleft(a\mright)\mright\}\mleft(n\mright) = i \ \ \ratio\Leftrightarrow \ \ \exists Y \left(H\mleft(a_0, Y\mright) \land \left(\left(\varphi\mleft(a, n, Y\mright) \land i = 1\right) \lor \left(\neg \varphi\mleft(a, n, Y\mright) \land i = 0\right)\right)\right)\]
    is $\Sigma^1_1$, where $\varphi\mleft(a, n, Y\mright) = a_1 \in \Rec^Y \land n \in R^Y_{a_1}$.
\end{proof}

\subsection{Realisability of a global well-ordering}
\label{subsec:realisability-cm-gwo-in-atr}

The $\Sigma^1_1$ function $\mathrm{cmp}$ defined in \autoref{cor:well-order-cmp} is a well-founded preorder on $\mathcal{O}$. This easily induce a relation\footnote{We used the most natural definitions when defining $\mathrm{cmp}$ in \autoref{cor:well-order-cmp} and when defining $\mathrm{cmp}^*$ here. The consequence is that $\mathrm{cmp}$ is a non-strict (i.e.\ reflective) ordering on $\mathcal{O}$ but $\mathrm{cmp}^*$ is a strict (i.e.\ asymmetric) ordering on $\mathcal{H}$.} $\mathrm{cmp}^* \in \left[\mathcal{H}, S_2\right]$  such that $\mleft\{\mleft\{\mathrm{cmp}^*\mright\}\mleft(a\mright)\mright\}\mleft(b\mright) = 1$ if and only if both
\begin{enumerate}[label=(\roman*)]
    \item $\mleft\{\mleft\{\mathrm{cmp}\mright\}\mleft(a_0\mright)\mright\}\mleft(b_0\mright) = 1 \land \left(\mleft\{\mleft\{\mathrm{cmp}\mright\}\mleft(b_0\mright)\mright\}\mleft(a_0\mright) = 0 \lor \left(\mleft\{\mleft\{\mathrm{cmp}\mright\}\mleft(b_0\mright)\mright\}\mleft(a_0\mright) = 1 \land b < a\right)\right)$,
    \item $b \in \mathcal{H}$, i.e.\ $\exists Y \left(H\mleft(b_0, Y\mright) \land b_1 \in \Rec^Y\right)$.
\end{enumerate}
Observe that our definition of $\mathrm{cmp}$ ensures that $\mleft\{\mleft\{\mathrm{cmp}\mright\}\mleft(a_0\mright)\mright\}\mleft(b_0\mright) = 1$ implies $b_0 \in \mathcal{O}$ and hence the corresponding $\mathrm{H}$-set exists, thus the negation of the conjunction $\text{(i)} \land \text{(ii)}$ above is
\[\neg \text{(i)} \lor \left(\text{(i)} \land \exists Y \left(H\mleft(b_0, Y\mright) \land b_1 \not\in \Rec^Y\right)\right).\]
This is $\Sigma^1_1$ and we can set $\mleft\{\mleft\{\mathrm{cmp}^*\mright\}\mleft(a\mright)\mright\}\mleft(b\mright) = 0$ when this holds, concluding that the $\Sigma^1_1$ function $\mathrm{cmp}^*$ is well-defined.

\begin{proposition}[$\mathrm{ATR}_0$]
    \label{prop:hyperarithmetic-indices-wo}
    The relation $\prec$ on $\mathcal{H}$ given by
    \[b \prec a \ \ \ratio\Leftrightarrow \ \ \mleft\{\mleft\{\mathrm{cmp}^*\mright\}\mleft(a\mright)\mright\}\mleft(b\mright) = 1\]
    is a well-ordering.
\end{proposition}

\begin{proof}
    $\prec$ is a linear ordering simply because any two elements in $\mathcal{O}$ are comparable.

    Suppose it is not a well-ordering, i.e.\ we have a second-order object $f : \mathbb{N} \rightarrow S_2^*$ that is a strictly decreasing function. Then $f_0 : \mathbb{N} \rightarrow \mathcal{O}$ given by $f_0\mleft(0\mright) = f\mleft(n\mright)_0$ satisfies
    \[\otyp\mleft(f_0\mleft(0\mright)\mright) \geq \otyp\mleft(f_0\mleft(1\mright)\mright) \geq \cdots.\]

    In other words, we have $\otyp\mleft(f_0\mleft(0\mright)\mright) \geq \otyp\mleft(f_0\mleft(n\mright)\mright)$ for any $n \in \mathbb{N}$, i.e.\ there exists a unique $\gamma_n \leq_\mathcal{O} f_0\mleft(0\mright)$ such that $\otyp\mleft(f_0\mleft(n\mright)\mright) = \otyp\mleft(\gamma_n\mright)$. Now, fix $\Omega \subseteq \mathbb{N}$ such that $H\mleft(f_0\mleft(0\mright), \Omega\mright)$, then all the comparison statements will be arithmetic with parameter $\Omega$. Thus, by arithmetic comprehension, there exists a function $f' : n \mapsto \gamma_n$ and we immediately have $f_0\mleft(0\mright) = f'\mleft(0\mright) \geq_\mathcal{O} f'\mleft(1\mright) \geq_\mathcal{O} \cdots$. Since $f_0\mleft(0\mright) \in \mathcal{O}$, $f'$ must stabilise after some $N$, implying that $\otyp\mleft(f_0\mleft(m\mright)\mright) = \otyp\mleft(f_0\mleft(n\mright)\mright)$ for any $m, n > N$.

    Finally, for any $m > N$, we have $f\mleft(m\mright) \succ f\mleft(m + 1\mright)$ yet $\otyp\mleft(f\mleft(m\mright)_0\mright) = \otyp\mleft(f\mleft(m + 1\mright)_0\mright)$, thus we must have $f\mleft(m\mright) > f\mleft(m + 1\mright)$ by the definition of $\mathrm{cmp}^*$. Thus, $f$ restricted to $\left\{m : m > N\right\}$ is a strictly decreasing sequence of natural numbers. This contradicts $\Sigma^0_0$-induction. Therefore, such $f$ cannot exist and $\prec$ is indeed a well-ordering.
\end{proof}

With the set-up so far, we will show that $\mathrm{ATR}$ can realise the axioms $\left(\mathrm{W}1\right), \left(\mathrm{W}3\right)$ in $\mathrm{GWO}$. As later parts of this section shall show, the entire extension $\mathrm{CM} + \mathrm{GWO}$ is unfortunately too strong and has to be dealt with in the system $\mathrm{BI}$ instead.

To interpret axioms $\left(\mathrm{W}1\right), \left(\mathrm{W}3\right)$ in $\mathrm{ATR}$, we shall use the same realisability conditions as in \autoref{subsec:realisability-cm}, but fix the interpretation of $\prec$ as
\[\mleft\{r\mright\}\mleft(e\mright) = \mleft\{\mleft\{\mathrm{cmp}^*\mright\}\mleft(\mleft\{\mathrm{hyp}\mright\}\mleft(e_1\mright)\mright)\mright\}\mleft(\mleft\{\mathrm{hyp}\mright\}\mleft(e_0\mright)\mright).\]

As stipulated in \autoref{thm:conv-s2-to-hyp}, this definition indeed satisfies the requirement above that
\[\forall e, e', f, f' \in S_2 \left(e \equiv f \land e' \equiv f' \rightarrow \mleft\{r\mright\}\mleft(\tuple{e, e'}\mright) = \mleft\{r\mright\}\mleft(\tuple{f, f'}\mright)\right).\]

Now, first observe that the axiom set $\left(\mathrm{W}1\right)$ stating that $\prec$ is a linear ordering is simply interpreted as propositional compositions of atomic formulae of the form $\mleft\{\mathrm{hyp}\mright\}\mleft(e\mright) \prec \mleft\{\mathrm{hyp}\mright\}\mleft(f\mright)$. Since this is indeed a linear ordering modulo extensional equality on $S_2$, these axioms are realised by any $\Sigma^1_1$ function with appropriate domain and range. We also have

\begin{proposition}
    The axiom of countable initial segments $\left(\mathrm{W}3\right)$ is realised:
    \[\forall X \, \exists Z \, \forall Y \left(Y \prec X \rightarrow \exists n \, Y = \left(Z\right)_n\right).\]
\end{proposition}

\begin{proof}
    Consider any $x \in S_2$. We make use of the fact that $\mleft\{\mathrm{cmp}^*\mright\}\mleft(\mleft\{\mathrm{hyp}\mright\}\mleft(x\mright)\mright) \in S_2$ is a $\Sigma^1_1$ function defining the initial segment $W = \left\{y \in \mathcal{H} : y \prec \mleft\{\mathrm{hyp}\mright\}\mleft(x\mright)\right\}$. Therefore, the desired $Z$ can simply be constructed as the code $z \in S_2$ given by
    \[\mleft\{z\mright\}\mleft(\tuple{m, n}\mright) = \left\{\begin{aligned}
             & \mleft\{\mleft\{\mathrm{inv}\mright\}\mleft(m\mright)\mright\}\mleft(n\mright) &  & \text{if} \ \mleft\{\mleft\{\mathrm{cmp}^*\mright\}\mleft(\mleft\{\mathrm{hyp}\mright\}\mleft(x\mright)\mright)\mright\}\mleft(m\mright) = 1, \\
             & 0                                                                              &  & \text{if} \ \mleft\{\mleft\{\mathrm{cmp}^*\mright\}\mleft(\mleft\{\mathrm{hyp}\mright\}\mleft(x\mright)\mright)\mright\}\mleft(m\mright) = 0.
        \end{aligned}\right.\]

    Our construction ensures that if any $y \in S_2$ satisfies $\mleft\{\mathrm{hyp}\mright\}\mleft(y\mright) \prec \mleft\{\mathrm{hyp}\mright\}\mleft(x\mright)$, then $\left(Z\right)_{\mleft\{\mathrm{hyp}\mright\}\mleft(y\mright)} = Y$. Therefore, the $\forall Y \left(Y \prec X \rightarrow \exists n \, Y = \left(Z\right)_n\right)$ is simply realised by
    \[\lambda y. \lambda d. \tuple{\mleft\{\mathrm{hyp}\mright\}\mleft(y\mright), *},\]
    where $*$ is any function with appropriate domain and range.
\end{proof}

\subsection{Unrestricted transfinite induction}
\label{subsec:unrestricted-ti}

Finally, we work on the axiom scheme of unrestricted transfinite induction $\left(\mathrm{W}2\right)$:
\[\forall X \left(\forall Y \left(Y \prec X \rightarrow \varphi\mleft(Y\mright)\right) \rightarrow \varphi\mleft(X\mright)\right) \rightarrow \forall X \, \varphi\mleft(X\mright)\]
for any third-order formula $\varphi\mleft(X\mright)$.

To realise this, we move to the stronger second-order theory $\mathrm{BI}$, that is, the axioms in $\mathrm{ACA}_0$ plus the scheme of transfinite induction along any countable well-ordering:
\[\forall X \, \forall {<_X} \left(\WO\mleft(X, {<_X}\mright) \rightarrow \forall x \in X \left(\forall y <_X x \ \varphi\mleft(y\mright) \rightarrow \varphi\mleft(x\mright)\right) \rightarrow \forall x \in X \, \varphi\mleft(x\mright)\right)\]
for any formula $\varphi\mleft(x\mright)$, where $\WO\mleft(X, {<_X}\mright)$ denotes the $\Pi^1_1$-formula asserting that $<_X$ is a linear ordering on $X$ and there exists no strictly decreasing function $f : \mathbb{N} \rightarrow X$. This is a further extension of the classical theories we worked with above due to corollary VII.2.19 in Simpson \cite{simpson09-soa}:

\begin{fact*}
    $\mathrm{BI}$ implies $\mathrm{ATR}$.
\end{fact*}

By working in the theory $\mathrm{BI}$, we additionally know how to iterate $\Sigma^1_1$ functions along any well-ordering:

\begin{theorem}[$\mathrm{BI}$]
    Given $x, y \in S_2$ and write
    \[X = \left\{n : \mleft\{x\mright\}\mleft(n\mright) = 1\right\}, \qquad {<_X} = \left\{n : \mleft\{y\mright\}\mleft(n\mright) = 1\right\},\]
    which exist as hyperarithmetic sets. Suppose that $\WO\mleft(X, {<_X}\mright)$. Fix some formula $\varphi\mleft(a, z\mright)$ and suppose $f$ is a $\Sigma^1_1$ function code such that for any $a \in X$ and $e \in \mathbb{N}$, if for any $b \in X$, $b <_X a$, we have $b \in \dom\mleft(e\mright)$ and $\varphi\mleft(b, \mleft\{e\mright\}\mleft(b\mright)\mright)$, then $\tuple{a, e} \in \dom\mleft(f\mright)$ and $\varphi\mleft(a, \mleft\{f\mright\}\mleft(\tuple{a, e}\mright)\mright)$. Then there is a $\Sigma^1_1$ function code $g$ such that
    \begin{itemize}
        \item $X \subseteq \dom\mleft(g\mright)$,
        \item $\varphi\mleft(a, \mleft\{g\mright\}\mleft(a\mright)\mright)$ for any $a \in X$,
        \item $\mleft\{g\mright\}\mleft(a\mright) = \mleft\{f\mright\}\mleft(\tuple{a, g}\mright)$ for any $a \in X$.
    \end{itemize}
    Additionally, $g$ can be computed from $f$, $x$ and $y$ uniformly by a $\Sigma^1_1$ coded function.
\end{theorem}

\begin{proof}
    Let $d$ be a $\Sigma^1_1$ code that computes $\lambda h. \lambda a. \mleft\{h\mright\}\mleft(\tuple{a, h}\mright)$. We now simply compute a code for the following formula
    \begin{align*}
        \mleft\{h\mright\}\mleft(a\mright) = z \ \ \ratio\Leftrightarrow & \ \ \mleft\{x\mright\}\mleft(a_0\mright) = 1                                                                                                                                                                                         \\
                                                                         & {} \land \forall b \left(\mleft\{y\mright\}\mleft(\tuple{b, a_0}\mright) = 1 \rightarrow \mleft\{a_1\mright\}\mleft(\tuple{b, a_1}\mright) = \mleft\{f\mright\}\mleft(\tuple{b, \mleft\{d\mright\}\mleft(a_1\mright)}\mright)\right) \\
                                                                         & {} \land z = \mleft\{f\mright\}\mleft(\tuple{a_0, \mleft\{d\mright\}\mleft(a_1\mright)}\mright),
    \end{align*}
    and we let $g = \mleft\{d\mright\}\mleft(h\mright)$. We have yet to show that the formula for $h$ on top defines a function indeed, but it will suffice if we can show the following formula
    \[\psi\mleft(a\mright) \coloneqq a \in \dom\mleft(g\mright) \land \mleft\{g\mright\}\mleft(a\mright) = \mleft\{f\mright\}\mleft(\tuple{a, g}\mright) \land \varphi\mleft(a, \mleft\{g\mright\}\mleft(a\mright)\mright)\]
    for all $a \in X$. Notice that suppose we have $\psi\mleft(b\mright)$ for all $b <_X a$, then for each $b$ we have $\tuple{b, h} \in \dom\mleft(h\mright)$ and indeed $\mleft\{h\mright\}\mleft(\tuple{b, h}\mright) = \mleft\{f\mright\}\mleft(\tuple{b, g}\mright)$ satisfies $\varphi\mleft(b, \mleft\{h\mright\}\mleft(\tuple{b, h}\mright)\mright)$. Thus $\tuple{a, g} \in \dom\mleft(f\mright)$ and $\mleft\{h\mright\}\mleft(\tuple{a, h}\mright) = z$ is satisfied by the unique $z$ such that $z = \mleft\{f\mright\}\mleft(\tuple{a, g}\mright)$, which also satisfies $\varphi\mleft(a, z\mright)$ by our assumption. Therefore, using $\mathrm{BI}$ on $\psi$, we can immediately conclude $\psi\mleft(a\mright)$ for all $a \in X$.

    Finally, the fact that $g$ can be computed using a $\Sigma^1_1$ function code simply follows from \autoref{prop:sigma11-smn}, treating the inputs as parameters for the function $h$.
\end{proof}

Therefore, assuming $\mathrm{BI}$, we can finally show that the realisability conditions introduced above satisfy transfinite induction for any formula $\varphi\mleft(X\mright)$ on the well-ordering $\prec$:

\begin{lemma}[$\AC$]
    \label{lem:realise-equivalent-parameters}
    For any $d \in \mathbb{N}$ and $e, f \in S_2$ such that $e \equiv f$, if $d \Vdash \varphi\mleft(e\mright)$ for some formula $\varphi\mleft(X\mright)$, then likewise $d \Vdash \varphi\mleft(f\mright)$.

    Specifically, assuming $\mathrm{ATR}$, then $d \Vdash \varphi\mleft(e\mright)$ if and only if $d \Vdash \varphi\mleft(\mleft\{\mathrm{inv}\mright\}\mleft(\mleft\{\mathrm{hyp}\mright\}\mleft(e\mright)\mright)\mright)$.
\end{lemma}

\begin{proof}
    Observe that when $\varphi$ is atomic, $d$ is irrelevant and this holds simply because our realisability conditions realise the Leibniz law. Thus, this holds for any formula $\varphi\mleft(X\mright)$ by simple induction on the complexity of $\varphi$.
\end{proof}

\begin{proposition}[$\mathrm{BI}$]
    The axiom of unrestricted transfinite induction $\left(\mathrm{W}2\right)$ is realised:
    \[\forall X \left(\forall Y \left(Y \prec X \rightarrow \varphi\mleft(Y\mright)\right) \rightarrow \varphi\mleft(X\mright)\right) \rightarrow \forall X \, \varphi\mleft(X\mright).\]
\end{proposition}

\begin{proof}
    For any $f \Vdash \forall X \left(\forall Y \left(Y \prec X \rightarrow \varphi\mleft(Y\mright)\right) \rightarrow \varphi\mleft(X\mright)\right)$ and any $x \in S_2$, we need to compute some $d$ such that $d \Vdash \varphi\mleft(x\mright)$. Again use that $\mleft\{\mathrm{cmp}^*\mright\}\mleft(\mleft\{\mathrm{hyp}\mright\}\mleft(x\mright)\mright) \in S_2$ is a $\Sigma^1_1$ function and we can consider the set
    \[W = \left\{y \in \mathcal{H} : \mleft\{\mleft\{\mathrm{cmp}^*\mright\}\mleft(\mleft\{\mathrm{hyp}\mright\}\mleft(x\mright)\mright)\mright\}\mleft(y\mright) = 1 \land \mleft\{\mathrm{hyp}\mright\}\mleft(\mleft\{\mathrm{inv}\mright\}\mleft(y\mright)\mright) = y\right\}.\]
    Notice that for any $y \in S_2$, we have $\mleft\{\mathrm{hyp}\mright\}\mleft(\mleft\{\mathrm{inv}\mright\}\mleft(\mleft\{\mathrm{hyp}\mright\}\mleft(y\mright)\mright)\mright) = \mleft\{\mathrm{hyp}\mright\}\mleft(y\mright)$ by stipulation in \autoref{thm:conv-s2-to-hyp}, thus $\mleft\{\mathrm{hyp}\mright\}\mleft(y\mright) \in W$ as long as $\mleft\{\mathrm{hyp}\mright\}\mleft(y\mright) \prec \mleft\{\mathrm{hyp}\mright\}\mleft(x\mright)$.

    Here we have $\WO\mleft(W, \prec\mright)$ by \autoref{prop:hyperarithmetic-indices-wo}, thus we can iterate a $\Sigma^1_1$ function on $W$. We will iterate the function $g'$ given by
    \[\lambda a. \mleft\{\mleft\{f\mright\}\mleft(\mleft\{\mathrm{inv}\mright\}\mleft(a_0\mright)\mright)\mright\}\mleft(\mleft\{d\mright\}\mleft(a_1\mright)\mright)\]
    where $d$ computes $\lambda e. \lambda z. \lambda h. \mleft\{e\mright\}\mleft(\mleft\{\mathrm{hyp}\mright\}\mleft(z\mright)\mright)$. Observe that for $y \in W$ and $e \in \mathbb{N}$, if for any $z \in W$, $z \prec y$, we have $\mleft\{e\mright\}\mleft(z\mright) \Vdash \varphi\mleft(\mleft\{\mathrm{inv}\mright\}\mleft(z\mright)\mright)$, then
    \[\mleft\{d\mright\}\mleft(e\mright) \Vdash \forall Z \left(Z \prec \mleft\{\mathrm{inv}\mright\}\mleft(y\mright) \rightarrow \varphi\mleft(Z\mright)\right)\]
    by \autoref{lem:realise-equivalent-parameters} and hence $\mleft\{g'\mright\}\mleft(\tuple{y, e}\mright) \Vdash \varphi\mleft(\mleft\{\mathrm{inv}\mright\}\mleft(y\mright)\mright)$. Therefore, iterating $g'$ on $W$ will produce a $\Sigma^1_1$ function code $g$ such that for any $y \in W$, $\mleft\{g\mright\}\mleft(y\mright) \Vdash \varphi\mleft(\mleft\{\mathrm{inv}\mright\}\mleft(y\mright)\mright)$. Then $\mleft\{g'\mright\}\mleft(\tuple{\mleft\{\mathrm{hyp}\mright\}\mleft(x\mright), g}\mright) \Vdash \varphi\mleft(x\mright)$ by \autoref{lem:realise-equivalent-parameters} again and is precise the number we need to compute.
\end{proof}

We have thus shown that in $\mathrm{BI}$, our realisability conditions interpret all intuitionistic consequences of $\mathrm{CM} + \mathrm{GWO}$.

\subsection{An ordinal analysis of \texorpdfstring{$\mathrm{CM} + \mathrm{GWO}$}{CM + GWO}}

We say that a formula $\varphi\mleft(X, x\mright)$ in the language of $\mathrm{CM}$ is a \emph{strictly positive first-order formula}, possibly with other free variables, if it is built up according the following rules:
\begin{itemize}
    \item $\varphi$ can contain any atomic formulae in the language of $\mathrm{PA}$, as well as atomic formulae of the form $t \in_1 X$ where $t$ is an arbitrary first-order term;
    \item $\varphi$ can contain the propositional constant ``bottom'': $\bot$;
    \item $\varphi$ can contain connectives $\land$, $\lor$, $\rightarrow$, with the restriction that $X$ cannot occur on the left of the connective $\rightarrow$;
    \item $\varphi$ can contain quantifiers $\forall y$ and $\exists y$ where $y$ is any first-order variable.
\end{itemize}

Let $\theta\mleft(x\mright)$ be any formula, then we let $\varphi\mleft(\theta, x\mright)$ denote the result of substituting every occurrence of $t \in_1 X$ inside $\varphi\mleft(X, x\mright)$ by $\theta\mleft(t\mright)$, with appropriate renaming so that any other free variable in $\theta\mleft(x\mright)$ is not bound by any quantifiers in $\varphi$.

We shall prove the following:

\begin{theorem}
    \label{thm:inductive-definition}
    For any strictly positive first-order formula $\varphi\mleft(X, x\mright)$, there exists $\theta_\varphi\mleft(x\mright)$ such that
    \[\mathrm{CM} + \mathrm{GWO} \vdash \forall x \left(\varphi\mleft(\theta_\varphi, x\mright) \rightarrow \theta_\varphi\mleft(x\mright)\right),\]
    and also for any formula $\eta\mleft(x\mright)$
    \[\mathrm{CM} + \mathrm{GWO} \vdash \forall x \left(\varphi\mleft(\eta, x\mright) \rightarrow \eta\mleft(x\mright)\right) \rightarrow \forall x \left(\theta_\varphi\mleft(x\mright) \rightarrow \eta\mleft(x\mright)\right).\]
\end{theorem}

This will involve a few lemmata:

\begin{lemma}[Positivity]
    \label{lem:pos}
    For any formulae $\theta\mleft(x\mright), \eta\mleft(x\mright)$,
    \[\mathrm{CM} \vdash \forall x \left(\theta\mleft(x\mright) \rightarrow \eta\mleft(x\mright)\right) \rightarrow \forall x \left(\varphi\mleft(\theta, x\mright) \rightarrow \varphi\mleft(\eta, x\mright)\right).\]
    This trivially implies $X \subseteq Y \rightarrow \forall x \left(\varphi\mleft(X, x\mright) \rightarrow \varphi\mleft(Y, x\mright)\right)$ as well.
\end{lemma}

\begin{proof}
    This follows from an easy induction on the complexity of $\varphi\mleft(X, x\mright)$.
\end{proof}

\begin{lemma}
    \label{lem:set-witnesses}
    For any formula $\theta\mleft(x\mright)$,
    \[\mathrm{CM} \vdash \forall x \left(\varphi\mleft(\theta, x\mright) \rightarrow \exists X \left(\forall z \in_1 X \, \theta\mleft(z\mright) \land \varphi\mleft(X, x\mright)\right)\right)\]
\end{lemma}

\begin{proof}
    We show this by induction on the complexity of $\varphi\mleft(X, x\mright)$:

    Firstly, if $\varphi\mleft(X, x\mright)$ is an atomic formula not of the form $t \in_1 X$, then $\varphi\mleft(\theta, x\mright) \rightarrow \varphi\mleft(\varnothing, x\mright)$ trivially; if $\varphi\mleft(X, x\mright)$ is $t \in_1 X$, then $\varphi\mleft(\theta, x\mright) \rightarrow \forall z \in_1 \left\{t\right\} \, \theta\mleft(z\mright) \land \varphi\mleft(\left\{t\right\}, x\mright)$.

    The case where $\varphi$ is $\bot$ is vacuously true.

    Now, suppose that $\varphi$ is of the form $\psi \land \chi$ and $\varphi\mleft(\theta, x\mright)$ holds, then by the inductive hypothesis we have $X$ and $Y$ such that
    \[\forall z \in_1 X \cup Y \, \theta\mleft(z\mright) \land \psi\mleft(X, x\mright) \land \chi\mleft(Y, x\mright).\]
    It follows from \autoref{lem:pos} that we also have $\psi\mleft(X \cup Y, x\mright) \land \chi\mleft(X \cup Y, x\mright)$.

    Suppose $\varphi$ is of the form $\psi \rightarrow \chi$ and $\varphi\mleft(\theta, x\mright)$ holds, we know from stipulation that $\psi\mleft(x\mright)$ must not contain $X$. By the inductive hypothesis, we have
    \[\psi\mleft(x\mright) \rightarrow \exists X \left(\forall z \in_1 X \, \theta\mleft(z\mright) \land \chi\mleft(X, x\mright)\right).\] Since $\psi\mleft(x\mright)$ is arithmetic, hence decidable, we can take $Y$ to be $X$ above if $\psi\mleft(x\mright)$ holds and $\varnothing$ otherwise. Then we always have $\forall z \in_1 Y \, \theta\mleft(z\mright)$, and we also have $\psi\mleft(x\mright) \rightarrow \chi\mleft(Y, x\mright)$.

    Suppose $\varphi$ is of the form $\forall y \, \psi$ and $\varphi\mleft(\theta, x\mright)$ holds, then by the inductive hypothesis
    \[\forall y \, \exists X \left(\forall z \in_1 X \, \theta\mleft(z\mright) \land \psi\mleft(X, x, y\mright)\right).\]
    By dependent choice, we can find $Z$ such that each $\left(Z\right)_y$ satisfies $\forall z \in_1 \left(Z\right)_y \, \theta\mleft(z\mright) \land \psi\mleft(\left(Z\right)_y, x, y\mright)$ and take $Y = \bigcup_y \left(Z\right)_y$. Then $\forall z \in_1 Y \, \theta\mleft(z\mright)$ and by \autoref{lem:pos}, we have $\forall y \, \psi\mleft(Y, x, y\mright)$ as desired.

    Finally, the cases where $\varphi$ is of the form $\psi \lor \chi$ or $\exists y \, \psi$ are trivial.
\end{proof}

\begin{proof}[Proof of \autoref{thm:inductive-definition}]
    Let $\psi\mleft(X, Z, S\mright)$ assert that both of the following holds:
    \begin{itemize}
        \item $\forall Y \preceq X \, \exists n \, Y = \left(Z\right)_n$;
        \item $\forall n \left(\left(Z\right)_n \preceq X \rightarrow \forall x \left(x \in_1 \left(S\right)_n \leftrightarrow \varphi\mleft(\bigcup \left\{\left(S\right)_m : \left(Z\right)_m \prec \left(Z\right)_n\right\}, x\mright)\right)\right)$.
    \end{itemize}
    By transfinite induction on all second-order objects, it is easy to show that for each $X$, a pair $Z, S$ exists and is unique in the following sense:
    \[\psi\mleft(X, Z, S\mright) \land \psi\mleft(X, Z', S'\mright) \rightarrow \forall n, m \left(\left(Z\right)_n \preceq X \land \left(Z\right)_n = \left(Z'\right)_m \rightarrow \left(S\right)_n = \left(S'\right)_m\right).\]

    Now, we let $\theta_\varphi\mleft(x\mright) \coloneqq \exists X, Z, S \left(\psi\mleft(X, Z, S\mright) \land \exists n \left(X = \left(Z\right)_n \land x \in_1 \left(S\right)_n\right)\right)$. We will verify that $\mathrm{CM} + \mathrm{GWO}$ proves the desired properties of $\theta_\varphi\mleft(x\mright)$:

    Firstly, suppose that $\varphi\mleft(\theta_\varphi, x\mright)$ holds for some $x$. By \autoref{lem:set-witnesses}, we can find $A$ satisfying
    \[\forall x \in_1 A \, \theta_\varphi\mleft(x\mright) \land \varphi\mleft(A, x\mright).\]
    Using dependent choice, the former conjunct will implies the existence of $X, Z, S$ such that
    \[\forall x \in_1 A \left(\psi\mleft(\left(X\right)_x, \left(Z\right)_x, \left(S\right)_x\mright) \land \exists n \left(\left(X\right)_x = \left(\left(Z\right)_x\right)_n \land x \in_1 \left(\left(S\right)_x\right)_n\right)\right).\]
    Using the fact that the $\prec$-initial segment before each $\left(X\right)_x$ is countable and dependent choice again, we can easily find $Z$ such that
    \[\forall x \in_1 A \, \forall Y \left(Y \preceq \left(X\right)_x \rightarrow \exists n \, Y = \left(Z\right)_n\right),\]
    so taking $Y = \left\{n : \neg n \in_1 \left(Z\right)_n\right\}$ and we must have $\forall x \in_1 A \left(X\right)_x \prec Y$. Let $Z^*, S^*$ be such that $\psi\mleft(Y, Z^*, S^*\mright)$ holds, then it follows that
    \[A \subseteq \bigcup \left\{\left(S^*\right)_m : \left(Z^*\right)_m \prec Y\right\}.\]
    Immediately $\varphi\mleft(A, x\mright)$ and \autoref{lem:pos} will imply $\theta_\varphi\mleft(x\mright)$.

    Secondly, suppose that for some formula $\eta\mleft(x\mright)$, $\forall x \left(\varphi\mleft(\eta, x\mright) \rightarrow \eta\mleft(x\mright)\right)$ holds. In order to have $\forall x \left(\theta_\varphi\mleft(x\mright) \rightarrow \eta\mleft(x\mright)\right)$, it suffices to show by transfinite induction that
    \[\forall X, Z, S \left(\psi\mleft(X, Z, S\mright) \rightarrow \forall n \left(X = \left(Z\right)_n \rightarrow \forall x \in_1 \left(S\right)_n \eta\mleft(x\mright)\right)\right).\]
    This is essentially because, when $X = \left(Z\right)_n$, the inductive hypothesis will imply
    \[\forall x \in_1 \bigcup \left\{\left(S\right)_m : \left(Z\right)_m \prec X\right\} \eta\mleft(x\mright),\]
    thus for any $x \in_1 \left(S\right)_n$, i.e.\ any $x$ satisfying $\varphi\mleft(\bigcup \left\{\left(S\right)_m : \left(Z\right)_m \prec X\right\}, x\mright)$, \autoref{lem:pos} would imply $\varphi\mleft(\eta, x\mright)$ and thus $\eta\mleft(x\mright)$ as required.
\end{proof}

Since $\mathrm{HA}$ is trivially a subtheory of $\mathrm{CM}$, the theorem implies that $\mathrm{CM} + \mathrm{GWO}$ interprets the intuitionistic\footnote{Even though $\mathrm{CM} \vdash \mathrm{PA}$, we can only interpret the intuitionistic part here because we are translating the first-order predicate symbol $I_\varphi$ for the least fixed point of $\varphi$ into a non-decidable formula $\theta_\varphi\mleft(x\mright)$.} first-order theory of inductive definitions $\mathrm{ID}_1^i$. The proof-theoretic ordinal of $\mathrm{ID}_1^i$ is computed in \cite{buchholz-pohlers78-iterated-id} as $\left|\mathrm{ID}_1^i\right| = \theta_{\varepsilon_{\Omega + 1}}\mleft(0\mright) = \left|\mathrm{BI}\right|$. Since $\mathrm{CM} + \mathrm{GWO}$ can be realised in $\mathrm{BI}$, it is an immediate corollary that

\begin{corollary}
    $\mathrm{CM} + \mathrm{GWO} \equiv \mathrm{BI}$. In particular, we have $\left|\mathrm{CM} + \mathrm{GWO}\right| = \left|\mathrm{BI}\right| = \theta_{\varepsilon_{\Omega + 1}}\mleft(0\mright)$.
\end{corollary}

\section{Indispensability of countable initial segments}
\label{sec:countable-initial-segments}

In this section, we will point out an interesting contrast between the extension $\mathrm{CM} + \mathrm{GWO}$ and the more conventional attempts of attaching a global well-ordering to a base theory in classical logic. For example, over Kripke--Platek set theory, one can take G\"odel's constructible universe $L$ as a witness of that $\mathrm{KP}$ is equiconsistent with ``$\mathrm{KP} + \text{there is a global well-ordering of all sets}$''. As an immediate by-product, the first uncountable ordinal $\omega_1$, if it exists, is regular and hence an admissible ordinal in $\mathrm{KP}$, thus $L_{\omega_1}$ will be a model of
\[\mathrm{KP} + \text{there is a global well-ordering whose every initial segment is countable},\]
rendering this theory also equiconsistent with $\mathrm{KP}$.

Through the ordinal analysis in the sections above, we have shown that adding the global well-ordering axioms $\mathrm{GWO}$ to the intuitionistic system $\mathrm{CM}$, instead, results in an increase in the consistency strength. If we likewise consider the fragment $\mathrm{CM} + \left(\mathrm{W}1\right) + \left(\mathrm{W}2\right)$ where the axiom of countable initial segments $\left(\mathrm{W}3\right)$ is dropped, then one may suspect that the simpler truncation argument of $L_{\omega_1}$ above may be replicated in the intuitionistic context whereas the more complicated construction of a well-orderable ``constructible class'' naturally fails, so a quick guess would be that $\mathrm{CM} + \left(\mathrm{W}1\right) + \left(\mathrm{W}2\right)$ inherits the same consistency strength of the full theory $\mathrm{CM} + \mathrm{GWO}$, rendering the extra axiom $\left(\mathrm{W}3\right)$ redundant.

However, we will show that this is not the case, and the axiom $\left(\mathrm{W}3\right)$ is indeed indispensable in order to achieve the proof-theoretic ordinal $\theta_{\varepsilon_{\Omega + 1}}\mleft(0\mright)$. We shall do so by providing an interpretation of the extension $\mathrm{CM} + \left(\mathrm{W}1\right) + \left(\mathrm{W}2\right)$ in the classical theory $\mathrm{ATR}$.

Let us again consider the realisability conditions given in \autoref{subsec:realisability-cm}, where the interpretation of the well-ordering $\prec$ is instead fixed as
\[\mleft\{r'\mright\}\mleft(\tuple{e, f}\mright) = \left\{\begin{aligned}
         & 0 &  & \text{if} \ \mleft\{\mathrm{hyp}\mright\}\mleft(e\mright) < \mleft\{\mathrm{hyp}\mright\}\mleft(e\mright), \\
         & 1 &  & \text{otherwise}.
    \end{aligned}\right.\]
Recall from \autoref{thm:conv-s2-to-hyp} that for $e, f \in S_2$, $\mleft\{\mathrm{hyp}\mright\}\mleft(e\mright) = \mleft\{\mathrm{hyp}\mright\}\mleft(e\mright)$ if and only if $e \equiv f$, thus $r'$ satisfies the extensionality requirement that
\[\forall e, e', f, f' \in S_2 \left(e \equiv f \land e' \equiv f' \rightarrow \mleft\{r\mright\}\mleft(\tuple{e, e'}\mright) = \mleft\{r\mright\}\mleft(\tuple{f, f'}\mright)\right).\]

Once again, we see that the axiom set $\left(\mathrm{W}1\right)$ stating that $\prec$ is a linear ordering is interpreted as propositional compositions of numerical inequalities in the base theory, thus realised by any $\Sigma^1_1$ function with appropriate domain and range. We also have

\begin{proposition}[$\mathrm{ATR}$]
    Using this simple interpretation $r'$, the axiom of unrestricted transfinite induction $\left(\mathrm{W}2\right)$ is realised:
    \[\forall X \left(\forall Y \left(Y \prec X \rightarrow \varphi\mleft(Y\mright)\right) \rightarrow \varphi\mleft(X\mright)\right) \rightarrow \forall X \, \varphi\mleft(X\mright).\]
\end{proposition}

\begin{proof}
    Given any $f \Vdash \forall X \left(\forall Y \left(Y \prec X \rightarrow \varphi\mleft(Y\mright)\right) \rightarrow \varphi\mleft(X\mright)\right)$, the traditional recursion theorem for partial combinatory algebras allows us to compute some $\Sigma^1_1$ partial function $g$ such that
    \[\mleft\{g\mright\}\mleft(n\mright) = \mleft\{\mleft\{f\mright\}\mleft(\mleft\{\mathrm{inv}\mright\}\mleft(n\mright)\mright)\mright\}\mleft(\lambda y. \lambda e. \mleft\{g\mright\}\mleft(\mleft\{\mathrm{hyp}\mright\}\mleft(y\mright)\mright)\mright).\]

    We can show by induction on $\mathbb{N}$ that the following holds for every $x \in S_2$:
    \[\forall y \in S_2 \left(\mleft\{\mathrm{hyp}\mright\}\mleft(y\mright) \leq \mleft\{\mathrm{hyp}\mright\}\mleft(x\mright) \rightarrow \mleft\{\mathrm{hyp}\mright\}\mleft(y\mright) \in \dom\mleft(g\mright) \land \mleft\{g\mright\}\mleft(\mleft\{\mathrm{hyp}\mright\}\mleft(y\mright)\mright) \Vdash \varphi\mleft(y\mright)\right).\]
    Here, the inductive hypothesis implies
    \[\lambda y. \lambda e. \mleft\{g\mright\}\mleft(\mleft\{\mathrm{hyp}\mright\}\mleft(y\mright)\mright) \Vdash \forall Y \left(Y \prec x \rightarrow \varphi\mleft(Y\mright)\right),\]
    thus the construction of $g$ ensures that $\mleft\{\mathrm{hyp}\mright\}\mleft(x\mright) \in \dom\mleft(g\mright)$ and, by \autoref{lem:realise-equivalent-parameters}, $\mleft\{g\mright\}\mleft(\mleft\{\mathrm{hyp}\mright\}\mleft(x\mright)\mright) \Vdash \varphi\mleft(x\mright)$ as well. Hence, by induction, this means that
    \[\lambda x. \mleft\{g\mright\}\mleft(\mleft\{\mathrm{hyp}\mright\}\mleft(x\mright)\mright) \Vdash \forall X \, \varphi\mleft(X\mright). \qedhere\]
\end{proof}

Compare this with the results we obtained in \autoref{subsec:realisability-cm-gwo-in-atr}, we observe that both the fragments $\mathrm{CM} + \left(\mathrm{W}1\right) + \left(\mathrm{W}2\right)$ and $\mathrm{CM} + \left(\mathrm{W}1\right) + \left(\mathrm{W}3\right)$ are realisable in $\mathrm{ATR}$ only --- although with different realisers for the well-ordering $\prec$ --- but the full extension $\mathrm{CM} + \mathrm{GWO}$ requires $\mathrm{BI}$. Therefore, indeed all axioms in our extension are necessary for its consistency strength.

\section{A concluding remark on impredicativity}
\label{sec:impredicativity}

In this last section, we will return to focus on the transfinite induction scheme $\left(\mathrm{W}2\right)$,
\[\forall X \left(\forall Y \left(Y \prec X \rightarrow \varphi\mleft(Y\mright)\right) \rightarrow \varphi\mleft(X\mright)\right) \rightarrow \forall X \, \varphi\mleft(X\mright).\]

As we mentioned when defining the axioms, Weaver did not want this scheme to apply to an arbitrary third-order formula $\varphi\mleft(X\mright)$, since it renders the resulting system impredicative. In his own words in \cite{weaver09-cm},
\begin{quote}
    [transfinite induction] should only be asserted for formulas that do not contain $\prec$, for reasons having to do with the circularity involved in making sense of a relation that is well-ordered with respect to properties that are defined in terms of that relation.
\end{quote}

However, while the rejection of impredicative constructions is understandable from the philosophical stance Weaver chose to take, his suggestion quoted here fail to resolve the problem. If we read it in the usual sense, i.e.\ we forbid the primitive relation symbol $\prec$ from $\varphi$ but still allow the occurrences of any other free parameters, then given any formula $\varphi\mleft(X\mright)$, one can simply apply the transfinite induction scheme to the formula $\tilde{\varphi}\mleft(X\mright)$ instead, where every occurrence of $X \prec Y$ has been replaced by $\tuple{X, Y} \in \mathbf{W}$, given $\mathbf{W}$ being an unused, free third-order parameter. By decidable comprehension, we can later set
\[\mathbf{W} \ \ \coloneqq \ \ \left\{\tuple{X, Y} : X \prec Y\right\}\]
to derive transfinite induction for the original formula $\varphi$. Therefore, such superficial restriction on $\varphi$ does not actually restrict the theory in any ways.

However, if we are to read this proposal in a harsher way, by either forbidding the occurrence of any free (third-order) parameter in $\varphi$ at all, or maintaining that the axiom scheme of decidable comprehension is postulated for the restricted language of $\mathrm{CM}$ only and cannot be applied to formulae containing the new symbol $\prec$, then the utility of the extension will be severely reduced. For example, as we have explored in \autoref{sec:axiom-gwo}, the proof of traditional mathematical results like the existence of bases using the global well-ordering axioms will very often involve recursive constructions along the well-ordering and then applications of the induction axiom to prove properties of such constructions. It is unclear how any of those results, including the general version of Hahn--Banach theorem that Weaver himself claimed to have proved, can be established without accessing the initial segment, i.e.\ second-order sets $Y \prec X$ for a given $X$ in a meaningful way.

Moving away from Weaver's original proposal, it still remains open though whether one can find any interesting theory fragment lying between $\mathrm{CM}$ and $\mathrm{CM} + \mathrm{GWO}$ that both
\begin{itemize}
    \item has significantly weaker proof-theoretic strength than full $\mathrm{CM} + \mathrm{GWO}$; and
    \item can still prove most of the mathematical results we looked at in \autoref{sec:axiom-gwo}.
\end{itemize}

For example, if we restrict the transfinite induction scheme $\left(\mathrm{W}2\right)$ to only where the formula $\varphi\mleft(X\mright)$ is decidable, then this is equivalent to the single axiom
\begin{enumerate}
    \item[$\left(\mathrm{W}2'\right)$] $\forall X \left(\forall Y \left(Y \prec X \rightarrow Y \in_2 \mathbf{X}\right) \rightarrow X \in_2 \mathbf{X}\right) \rightarrow \forall X \, X \in_2 \mathbf{X}$
\end{enumerate}
by taking the third-order parameter $\mathbf{X} \coloneqq \left\{X : \varphi\mleft(X\mright)\right\}$. We can prove that

\begin{proposition}[$\mathrm{ATR}$]
    \label{prop:atr-realise-decidable-ti}
    The axiom of transfinite induction $\left(\mathrm{W}2'\right)$ is realised (using the realiser $r$ defined in \autoref{subsec:realisability-cm-gwo-in-atr} for the symbol $\prec$).
\end{proposition}

\begin{proof}
    Fix some $h \in S_3$ and suppose that $\forall X \left(\forall Y \left(Y \prec X \rightarrow Y \in_2 h\right) \rightarrow X \in_2 h\right)$ is realised. This implies
    \[\forall e \in S_2 \left(\forall f \in S_2 \left(\mleft\{\mathrm{hyp}\mright\}\mleft(f\mright) \prec \mleft\{\mathrm{hyp}\mright\}\mleft(e\mright) \rightarrow \mleft\{h\mright\}\mleft(f\mright) = 1\right) \rightarrow \mleft\{h\mright\}\mleft(e\mright) = 1\right),\]
    and it suffices to show that $\forall e \in S_2 \, \mleft\{h\mright\}\mleft(e\mright) = 1$ as then any $\Sigma^1_1$ function with appropriate domain and range would realise $\forall X \, X \in_2 h$.

    Suppose that this is not the case, then there exists some $e \in S_2$ such that $\mleft\{h\mright\}\mleft(e\mright) = 0$, and it follows that we have some $f \prec \mleft\{\mathrm{hyp}\mright\}\mleft(e\mright)$ satisfying $\mleft\{h\mright\}\mleft(\mleft\{\mathrm{inv}\mright\}\mleft(f\mright)\mright) = 0$. Fix $e$ and consider the set
    \[\left\{f \in \mathbb{N} : \mleft\{\mleft\{\mathrm{cmp}^*\mright\}\mleft(\mleft\{\mathrm{hyp}\mright\}\mleft(e\mright)\mright)\mright\}\mleft(f\mright) = 1 \land \mleft\{h\mright\}\mleft(\mleft\{\mathrm{inv}\mright\}\mleft(f\mright)\mright) = 0\right\}\]
    which exists by $\Delta^1_1$-comprehension since $\mleft\{\mathrm{cmp}^*\mright\}\mleft(\mleft\{\mathrm{hyp}\mright\}\mleft(e\mright)\mright)$ is a total $\Sigma^1_1$ function, then it must be non-empty and thus contains a $\prec$-minimal element. However, this contradicts the assumption above.
\end{proof}

Combined with results in \autoref{subsec:realisability-cm-gwo-in-atr}, this means that the fragment $\mathrm{CM} + \left(\mathrm{W}1\right) + \left(\mathrm{W}2'\right) + \left(\mathrm{W}3\right)$ is entirely realised $\mathrm{ATR}$. However, this fragment is still too weak, as one can observe that the crucial use of the transfinite induction axiom in \autoref{prop:zorns-lemma} requires one second-order existential quantifier in the formula $\varphi$, and there is no convenient prior justification that the specific $\varphi$ involved is indeed decidable. It shall be up to further work to determine whether the realisability result \autoref{prop:atr-realise-decidable-ti} can be extended to accommodate formulae of higher complexity.

\bibmain

\end{document}